\documentclass[12pt]{article}
\usepackage{amssymb}
\usepackage{amsmath}
\usepackage{hyperref}

\newcommand{\D}{\displaystyle}
\newcommand{\R}{\mathbb{R}}
\newcommand{\C}{\mathbb{C}}
\newcommand{\Z}{\mathbb{Z}}

\newtheorem{theorem}{Theorem}

\newtheorem{proposition}{Proposition}
\newtheorem{remark}{Remark}
\newtheorem{lemma}{Lemma}

\begin{document}

\title{SELFSIMILAR SOLUTIONS OF THE BINORMAL FLOW\\ AND THEIR STABILITY}
\author{Valeria Banica\footnote{}, Luis Vega\footnote{E-mail addresses: Valeria.Banica@univ-evry.fr, luis.vega@ehu.es }\footnote{First author was partially supported by the ANR project ``\'Etude qualitative des EDP dispersives". Second author was partially supported by the grant MTM 2007-62186 of MEC (Spain) and FEDER .}}


\vspace{4mm}
\begin{center}{{\large{\bf{SELFSIMILAR SOLUTIONS OF THE BINORMAL FLOW\\ AND THEIR STABILITY}}}}
\end{center}
\vspace{3mm}

\begin{center}
{\bf{Valeria Banica$^*$, Luis Vega$^\sharp$\footnote{First author was partially supported by the ANR projects ``\'Etude qualitative des EDP" and ``Equations de                       
Gross-Pitaevski, d'Euler, et ph\'enom\`enes de                                      
concentration". Second author was partially supported by the grant MTM 2007-62186 of MEC (Spain) and FEDER. E-mail addresses: Valeria.Banica@univ-evry.fr, mtpvegol@lg.ehu.es}\vspace{3mm}\\\tiny{$^*$D\'epartement de Math\'ematiques, Universit\'e d'Evry, France\\
$^\sharp$Departamento de Matematicas, Universidad del Pais Vasco, Spain }}}  
\end{center}

\abstract{We review some recent results concerning the evolution of a vortex filament and its relation to the cubic non-linear Schr\"odinger equation. Selfsimilar solutions and questions related to their stability are studied.\par Keywords: vortex filaments,  selfsimilarity, non-linear 
Schr\"odinger.}

\tableofcontents

\vspace{8mm}

\section{Introduction}

\subsection {The binormal flow}
This paper is mainly concerned with the so called Localized Induction approximation for the evolution of a vortex filament. In geometrical terms it is given by the binormal flow (BF): consider the 3d curve $\,\chi(s,t)\in\R^3\,$ with $\,t\in\R\,$ denoting time and $\,s\in\R\,$  the arclength parameter, then the flow is given by
  \begin{equation}
\label{I1}
\left\{\begin{array}{l}  \chi_t=\chi_s\,\land\,\chi_{ss}\\
  \chi(s,0)=\chi_0(s).
 \end{array}\right.
\end{equation}
Denoting $\,c\,$ the curvature and $\,b\,$ the binormal vector we get from the Frenet equations 
\begin{equation}
\label{I1bis}
\left\{\begin{array}{l}  \chi_t=cb\\
  \chi(s,0)=\chi_0(s).
 \end{array}\right.
\end{equation}

This flow appeared for the first time in a work in 1906 by L.S Da Rios  \cite {DaR} as a model which describes the flow of a vortex filament. That is to a say a vortex tube of infinitesimal cross section which propagates according to Euler equations. Let us call $\,\chi(s,t)\,$ a curve in $\,\R^3\,$ which can be either closed, and then we will consider $\,s\in[0,2\pi]\,$ and $\,2\pi$--periodic functions, or which extends to infinity in both directions. In that case $\,s\in\R\,$. Let us consider for simplicity of the exposition this second possibility.
  
  Assume that the vorticity $\,{w}\,$ of a given fluid of velocity $\,{u}\,$ is supported along the filament described by the curve$\,\chi(s,t)\,$. By that we mean $\,{w}\,$ is the singular vectorial measure
   $${w}=\Gamma Tds,$$
  where $\,\Gamma\,$ stands for the circulation,
  $$T=\chi_s,$$
  $\,s\,$ is arclength parametrization and $\,\Gamma\,$ is constant along the filament. 
    The velocity can be recovered from the vorticity through the Biot--Savart integral
  \begin{equation}{u}(P)=\D\frac{\Gamma}{4\pi}\D\int_{-\infty}^{\infty}\D\frac{\chi(s)-P}{|\chi(s)-P|^3}
\land T(s)ds.
\label{eq1}
\end{equation}


Typical examples of vortex structures are:
\begin{enumerate}
 \item
The straight line vortex,
$\chi(s)=(0,0,s)$. The velocity at a point $P=(x,y,z)$ with $z\neq0$  can be 
easily computed from  
(\ref{eq1})
$${u}(x,y,z)=\D\frac{\Gamma}{4\pi}(-y,x,0)\int_{-\infty}^{\infty}\D\frac{ds}{(x^2+y^2+(z-s)^2)^\frac32}=\D  
(-y,x,0)\frac{C}{x^2+y^2}.$$
  Hence the velocity around the filament is just a rotation of speed  a 
constant times the inverse of the distance to the filament. As a  consequence 
the filament remains 
stationary.
   \item The vortex ring, $\,\chi(s,0)=(\cos s,\sin s,0)$.
  \item The helical vortex, $\,\chi(s,0)=(\cos s/\sqrt2,\sin  
s/\sqrt2,/\sqrt2),$.-see \cite {LF} and 
\cite{Har}

  The computation of the velocity field outside of the filament is  much more 
delicate in the last two examples and in fact, it is better  to compute the so 
called vector potential using cylindric variables.  We refer the reader to  
\cite{AKO}, sections 2.5 and  2.6 for the  
details.
\end{enumerate}

Notice that  the Biot--Savart integral is perfectly defined for points $P$ outside of the filament. The idea is to compute $\,{u}(P)\,$ as $\,P\,$ approaches the filament to eventually identify that limit with the velocity of the filament. Also observe that the term inside of the integral in (\ref{eq1}) has enough decay at infinity. So in order to understand the leading behaviour  we can just  look at a neighbourhood of $\,\chi(s_0)$, the closest point in the filament to $\,P\,$. A quick look at the integral  is enough to show that the limit diverges as the inverse of the distance. In fact we see that the first approximation is
$$u(P)\sim\D\frac{\chi(s_0)-P}{|\chi(s_0)-P|}\land T(s_0).$$

This term behaves as in the case of example 1 of the straight filament. As we said the overall effect on the movement of the filament is null.  Doing the Taylor expansion around  that point, calling $\epsilon$ the distance of $P$ to the filament and removing the first term  we get (see for example \cite{Sa})
$$\begin{array}{rcl}
u(P)
&\sim&-c\Gamma\log\epsilon\,\chi_s(s_0)\land \chi_{ss}(s_0)+\cdots
\end{array}$$

Hence we see we have a logarithmic dependence on the distance to the filament on the direction of 
$$\chi_s\land\chi_{ss}=b\,(s_0).$$
After a renormalization of time we get
\begin{equation}
\label{eq2}
\chi_t=\chi_s\land\chi_{ss}.
\end{equation}

This computation is extremely crude and at best indicates that close to the filament the leading term is in the direction of the binormal. The general consense (see \cite {Sa}) about this approximation also known as local induction approximation, is that it is specially suitable for nummerical integration but it has serious flaws and is of limited physical interest. In any case it would be extremely useful to find a way of desingularizing the integral in (\ref{eq1}).

Our next step is to establish the connection between the binormal flow and non-linear Schr\"odinger equation. A first hint of such a connection comes by identifying  the normal plane to the curve $\,\chi(s)\,$ at a given point with the complex plane. Then if $n$ and $b$ are the normal and the binormal vector respectively we can write the binormal as
$$b=in.$$
Hence the binormal flow 
can be understood as
$$\D\frac 1i\,\chi_t=cn.$$
Finally if we remember that the ``curvature flow"  has as a linear approximation the heat flow, we can ``guess" that the underlying linear equation for (\ref{eq2}) will be Schr\"{o}dinger equation.
In fact something even stronger is true.

H. Hasimoto in \cite{Has} observed that if we define the complex valued function $\,u\,$ as
$$u(s,t)=c(s,t)\,e^{i\int_0^s\tau(s^{\prime},t)ds^{\prime}},$$
where $\,\tau\,$ stands for the torsion of the curve and $c$ for the the curvature, and if $\,\chi(s,t)\,$ solves (\ref{eq2}) then $\,u\,$ solves
\begin{equation}
\label{eq3}
iu_t+u_{ss}+\D\frac 12\left(|u|^2+A(t)\right)u=0,
\end{equation}
with
\begin{equation}\label{Hasimoto-coeff}
A(t)=\left(2\frac{c_{ss}-c\,\tau^2}{c}+c^2\right)(0,t).
\end{equation}

Conversely, let $A(t)$ be a time depending function, and $u$ a regular non vanishing solution of \eqref{eq3}. If we write
$$u(s,t)=c(s,t)e^{i\phi(s,t)},$$
by defining $\tau(s,t)=\phi_s(s,t)$, we have 
$$u(s,t)=c(s,t)e^{i\int_0^s\tau(s',t)ds' +\phi(0,t)}.$$
We obtain then a filament function 
$$\tilde{u}(s,t)=c(s,t)e^{i\int_0^s\tau(s',t)ds'},$$
satisfying \eqref{eq3} with $A(t)$ replaced by $A(t)+\phi_t(0,t)$. The real functions $c$ and $\tau$ verify then the intrinsic equations and the identity \eqref{Hasimoto-coeff} with $A(t)$ replaced by $A(t)+\phi_t(0,t)$. 

Notice that taking $\,\tilde u=u\exp\left(-\frac{i}{2}\int_0^t A(t^{\prime})dt^{\prime}\right)\,$ we get
\begin{equation}
\label{eq4}
i\tilde u_t+\tilde u_{ss}+\D\frac 12|\tilde u|^2\tilde u=0,
\end{equation}
which is the well known cubic non linear Schr\"odinger equation.

In doing this H. Hasimoto uses the Frenet frame $\,(T,n,b)\,$ and the corresponding Frenet equations
\begin{equation}
\label{Frenet}
\begin{array}{rclll}
T_s&=&&\,\,\,c\,n\,,\\
n_s&=&-c\,T&&+\,\tau \,b\,,\\
b_s&=&&-\,\tau\, n\,.
\end{array}
\end{equation}

There is no reason why to do this and more general frames $\,(T,e_1,e_2)\,$ can be used so that the assumption $\,c>0\,$ is not necessary, see for example \cite{Ko}. The final conclusion is still the same and equation (\ref{eq3}) is also obtained.

 Also it is important to observe that if in (\ref{eq2}) we differentiate with respect to $\,s\,$ in both sides we get, since $\chi_s=T$,
\begin{equation}
\label{eq5}
T_t=\chi_s\land \chi_{sss}=T\land T_{ss}.
\end{equation}

Hence
$$\D\frac d{dt}|\chi_s|^2=2\frac d{dt}\chi_s\cdot\chi_s=2\left(\chi_s\land\chi_{ss}\right)_s\cdot\chi_s=0,$$
and the arclength parametrization is preserved. This is generally admitted as one of the flaws of BF as an  approximation of (\ref{eq1}).

Notice that $\,|T|^2=1\,$ and $\,T\land T_{ss}=JD_sT_{s}\,$ where $\,J\,$ is the complex structure on the sphere and $D_s$ the covariant derivative. Therefore we can write (\ref{eq3}) as
$$T_t=JD_sT_{s}$$
which is nothing but the Schr\"odinger map onto the sphere in the one dimensional case. For higher dimensions we get
\begin{equation}
\label{eq6}
T_t=T\land \Delta T
\end{equation}
which also appeared in ferromagnetism where is referred as the Lipschitz--Landau equation.

Let us finish this short introduction to the binormal flow observing that Hasimoto transformation allows us to exhibit some particular examples of (\ref{eq2}):
\begin{enumerate}
  \item The straight line is just $\,c(s,t)=0,\quad\quad A(t)=0$.
  \item The circle is just $\,c(s,t)=1\quad\quad A(t)=-1$.
  \item The helix is just
  $$u(s,t)=e^{-itN^2} e^{iNs}\quad\quad A(t)=-1.$$
   \item The travelling wave solutions
   $$u(s,t)=e^{-itN^2} e^{iNs}Q(s-2Nt),\quad\quad A(t)=-1,$$
   with $Q(s)=\frac1{2\sqrt2}\, \sec h\,(s).$ This last example is particularly interesting because suggests the possibility of the existence of travelling wave solutions along          
vortex filaments. In fact a few years after H. Hasimoto's paper, E.J. Hopfinger and F.K. Browand \cite{HB} (see also \cite {MHR}) gave      
experimental evidence of such a possibility. 
\end{enumerate}

It is also interesting to recall some results by R. Klein, A. Majda and K. Damodaran \cite     
{KMD} about the interaction of more than one filament. They propose a mathematical model based on a system of coupling NLS. This model     
is obtained as an equilibrium between the interaction of a filament with itself, that is described just by the linear part of BF, and       
with the others, that is measured by a sum of culombic type potentials, see \cite{MB} for the details. This system allows for helixes as    
solutions. In the case of two or three helixes with the same vorticity and at the same distance the configuration is stable, see            
\cite{KPV1}. There is some experimental evidence that the case of four helixes could be unstable, see \cite{AKOS}.

Finally recall that  the $\,L^2$ norm is conserved in NLS. In our setting this becomes \begin{equation}
\label{eq7}
\D\int_{-\infty}^{\infty} c^2(s,t)\,ds=\int_{-\infty}^{\infty} |u(s,t)|^2\,ds,
\end{equation}
which can be understood as the kinetic energy. Notice that in the case of the helix we obtain a solution of infinite energy eventhough it has finite $\,L^2$--norm when consider as a periodic solution. As a filament one has to look at the complete helix because any vortex filament either ends in the boundary which confines the fluid or it has to go to infinity. This suggests as a natural question to consider solutions of (\ref{eq7}) with infinite energy. 
\subsection
{The mKdV flow}
The connection between NLS, which is a completely integrable system, and a geometric flow raised the natural question about the            
existence of more examples with similar connections. G. Lamb in the 70\'{}s (\cite{L}) was the first one in extending Hasimoto\'{}s idea, obtaining a remarkable generalization which include NLS and mKdV (even as complex valued functions) as possible limits. It is what is known as Hirota equation. Y. Fukumoto and T. Miyazaki \cite{FY} obtained the same flow when in da Rios approximation is included the existence of an axial flow along the filament. So it is still a 3d--flow.

Almost at the same time R.E. Goldstein and D.M. Petrich \cite{GoP1} proposed the mKdV hierarchy as the PDE of the curvature $\,k(s,t)\,$ of plane curve flows which preserve the enclosed area. In a later paper and in a complete formal fashion they relate it to the dynamics of a vortex patch in the plane. Recall that for plane fluids the vorticity is a scalar, see \cite{GoP2}. A vortex patch describes the situation of a compactly supported vorticity with a jump discontinuity through the boundary. Assume for simplicity that the vorticity is a constant $\,\Gamma\,$ in the patch $\,\Omega\,$ which is interior to the boundary $\,\partial\Omega\,$. Then N.J. Zabusky, M.H. Hughes and K.V. Roberts \cite{ZHR} obtain an equation for the boundary which is given by a singular integral coming from Biot-Savart law. After doing a formal Taylor expansion in this integral (see \cite{GoP2}) the first term is just a translation along the boundary. The second term is
$$z_t=-z_{sss}+\D\frac 32\overline{z}_s\,z_{ss}^2\qquad\qquad
z_t=k_s{n}-\D\frac 12k^2{T}.$$

If we call $\,k\,$ the curvature one gets
$$k_t=-k_{sss}-\D\frac 32k^2k_s$$
which is mKdV equation.

The biggest flaw of this model is that it is also arclength preserving  which 
is something that one cannot assume in the vortex--patch motion  due to the 
filamentation phenomena. In fact, since the early  simulations of Deem and 
Zabusky \cite{DZ}  there is a strong numerical  evidence (see \cite {Ar} and 
\cite{Bu2}) that small perturbations of a  circular patch can generate long and 
very thin fingers, and as a  consequence a big increasing of the total length 
of the boundary of  the patch is produced. As far as we know, nothing has been 
rigorously  proved so far that can support the connection of this geometric 
flow  with the one of the vortex patches, (see \cite{We} for a discussion about this point
and also about the connection of the travelling wave solutions of the  mKdV 
with the V-states of Deem and Zabusky introduced in  \cite{DZ}). However 
it would be interesting  to know if under some  assumptions, and possibly with 
some extra equations that keep track of  the evolution of the length of the 
curve, the mKdV can be useful to  describe some of the phenomena that appears 
in the dynamics of the  vortex 
patches.\newline

The rest of the paper will be devoted to the BF and its connection with NLS. 
In the next section we shall present some results on NLS equations with rough data, that will lead to some 
intuition on the existence of selfsimilar solutions for the binormal flow. Section \S 3 is devoted to the construction of
 some selfsimilar solutions which are examples of regular curves which develop a singularity in a later time. 
This singularity is either a corner or a 3d logarithmic spiral. Analogous examples but of plane curves can be given for the mKdV flow, at 
least with some smallness condition in the curvature. We refer the reader to the paper by G. Perelman and L. Vega [PeV] for  the 
corresponding results, and to the one by de la Hoz, \cite{Patxi}  for some 
numerical simulations concerning these 
examples. In the last section we shall study the stability of the selfsimilar solutions constructed in \S 3.1-\S 3.3.

\section{Cubic NLS with infinite energy}

\subsection{Some symmetries} 
 Recall the cubic NLS equation
\begin{equation}
\label{eq1.1}
\left\{\begin{array}{rcl}
iu_t+u_{ss}+\D\frac 12|u|^2u&=&0\\
u(s,0)&=&u_0(s).
\end{array}\right.
\end{equation}
This equation is scaling invariant:
\begin{itemize}
  \item If $\,\lambda>0\,$ then
  $\,u_{\lambda}=\lambda u(\lambda s,\lambda^2t)\,$ is a solution if so is $\,u$.
\end{itemize}
It is also galilean invariant:
\begin{itemize}
  \item If $\,N\in\R\,$ and $\,u\,$ is a solution then
  \begin{equation}
\label{eq1.2}
u_N(s,t)=e^{-itN^2+iNs}u(s-2Nt,t)
\end{equation}
is also a solution.
\end{itemize}

Notice that if $\widehat u$ denotes the Fourier transform then
$$\widehat{u_N}(\xi)=\widehat{u}(\xi-N).$$
Hence if $u$ belongs to the Sobolev space $H^k$,
$$\begin{array}{rcl}
\D\lim_{N\to\infty}\int\left| \widehat{u}_N\right|^2\left(1+|\xi|\right)^{2k}d\xi&=&
\D\lim_{N\to\infty}\int\left| \widehat u(\xi)\right|^2\left(1+|\xi-N|\right)^{2k}d\xi\\
&=&
\begin{cases}
      0& \text{if }k<0 \\
      \infty& \text{if }k>0.
\end{cases}
\end{array}$$

Assume that we have a result that assures that if $\,u_0\in H^k\,,$ then there exists a unique solution of (\ref{eq1.1}) up to time $\,t_0=t_0\left(\|u_0\|_{H^k}\right)\,$ with $\, t_0(r)\to\infty\,$ as $\,r\to 0\,$. Then if $\,k<0\,$ we will immediately obtain a global existence result from a local existence one just by using galilean transformations. This suggests to consider (\ref{eq1.1}) a supercritical equation from the point of view of galilean transformations and within the class of Sobolev spaces $\,H^k\,$ if $\,k<0$.

\subsection{NLS with rough data}
 In this section we want to review three results about NLS. 
 
 First of all we want to address the issue of selfsimilar solutions with respect to scaling within a general framework. 
 This was first studied by T. Cazenave and F.B. Weissler \cite{CW} and then extended by F. Planchon \cite{Pl1}, \cite{Pl2} 
 to NLS equations with nonlinearities $|u|^\alpha u$ which are $\,L^2\,$ critical or supercritical, 
  in the Besov spaces $\dot B^{\frac n2-\frac 2\alpha,\infty}_2$. But recall we are interested in 1D cubic NLS, which is 
  subrcritical in $\,L^2$. At this respect T. Cazenave, L. Vega and M.C. Vilela \cite{CVV} were able to 
  extend F. Planchon's result to nonlinearities which are subcritical in $\,L^2\,$ although the 1d cubic NLS is not included 
  in their result.

Let us consider the power like NLS
\begin{equation}
\label{eq3.1}
\left\{\begin{array}{rcl}
iu_t+\Delta u+\lambda|u|^\alpha u&=&0\qquad s\in\R^d\,,\quad t\in\R\,,\\\\
u(s,0)&=&u_0(s),
\end{array}\right.
\end{equation}
and focus on $\,d=1\,$ and $\,d=2\,$ and $\,\lambda=\pm1$. In \cite{CVV} the following result was proved. 

\begin{theorem}\label{3.T.1} Suppose
$$\D\frac{4(d+1)}{d(d+2)}<\alpha<\D\frac{4(d+1)}{d^2}.$$
There exists $\delta>0$ such that if $\,u_0\in\mathcal{S}^{\prime}(\R^d)\,$ satisfies $\,\|u_0\|_{\mathcal{Y}}<\delta\,$ then there exists a unique solution $\,u\in L^{\frac{\alpha(d+2)}2,\infty}\left(\R^{d+1}\right)\,$ of \eqref{eq3.1}  such that
$$\|u\|_{L^{\frac{\alpha(d+2)}2,\infty}(\R^{d+1})}\le 2\delta.$$

Here $\,\mathcal{Y}\,$ is defined as
$$\mathcal{Y}=\left\{\varphi\in\mathcal{S}^{\prime}(\R^d)\,;\,
e^{it\Delta}\varphi\in  L^{\frac{\alpha(d+2)}2,\infty}(\R^{d+1})\right\},$$
and \eqref{eq3.1} is understood as the associated integral equation:
$$u(t)=e^{it\Delta}u_0+i\lambda\int_0^te^{i(t-\tau)\Delta}|u|^\alpha u \,d\tau.$$
\end{theorem}

Notice that
\begin{itemize}
  \item For $d=1$, $\,\D\frac 83<\alpha<8\,$ and we would like to obtain $\alpha=2$.
  \item For $d=2$, $\,\D\frac 32<\alpha<3$.
  \end{itemize}

The key issue is as we see to choose a suitable space of initial data which is defined through the free evolution, namely
$$u_0\in \mathcal{Y}\quad\mbox{iff}\quad e^{it\Delta}u_0\in L^{p,\infty}(\mathbb{R}^{d+1}),$$
where $\,L^{p,\infty}(\mathbb{R}^{d+1})\,$ denotes the Lorentz space. These spaces are naturally defined by real interpolation which behaves poorly within the setting of mixed norm spaces $\,L^p(\mathbb{R},L^q(\mathbb{R}^d))\,$ and that is the reason for assuming $p=q$ in Theorem \ref{3.T.1}.

The other key ingredient is to use some Strichartz estimates for the operator
\begin{equation}
\label{eq3.2}
\mathcal{T}(F)=\int_0^t e^{i(t-\tau)\Delta} F(\cdot,\tau)d\tau
\end{equation}
within $\,\left(L^{p_1}, L^{p_2}\right)\,$ spaces which where unknown and involve pairs $\,(p_1,p_2)\,$ which are non admissible in the sense of J. Ginibre and G. Velo \cite{GiV}. The possibility of the existence of this type of inequalities was first observed by T. Cazenave and F.B. Weissler in \cite{CW}. Based on the work of M. Keel and T. Tao \cite{KeTa}, D. Foschi  \cite{F} and M.C. Vilela \cite{V} independently proved a general result at this respect although there are still open questions in this direction.

However the inequalities for the Duhamel operator are easier to prove if one relates $\,\mathcal{T}\,$ with Bochner--Riesz operators. These operators have been widely studied in Fourier analysis and are given in terms of the multiplier
$$m_{\gamma}(\xi)=\left(1-|\xi|^2\right)^{\gamma}_+,\quad \xi\in\R^d$$
with $\,(\cdot)_+\,$ denoting the positive part. Then the Bochner--Riesz operator of order $\,\gamma\,$ is defined as
$$BR_\gamma(f)^{\,\widehat{}\,}(\xi)=m_{\gamma}(\xi)\widehat{f}(\xi).$$

The index of interest is $\,\gamma=-1$, which is closely related to the more natural multiplier
\begin{equation}
\label{eq3.3}
H(\xi)=\left(1-|\xi|^2\right)^{-1},\quad \xi\in \R^d
\end{equation}
understood in the sense of the principal value, and which is linked to Helmholtz equation
\begin{equation}                                                                                                                            
\label{eq3.4}                                                                                                                               
\Delta u(x)+u(x)=f(x),\quad x\in \R^d.                                                                                                      
\end{equation}                                                                                                                              
The  connection with \eqref{eq3.2}  is clear because $\,\mathcal{T}F\,$ is a solution of
$$i\partial_tu+\Delta u=F,$$
and therefore, at least formally
$$\widehat{u}=\D\frac{1}{-\tau-|\xi|^2}\widehat{F}.$$

By scaling considerations one can always consider that $\,\mbox{supp}\,\widehat{F}(\xi,\tau)\subset \{|\xi|\le 1\}\,$ and then one mimics the available proof of the corresponding inequalities for solutions of \eqref{eq3.4} using the multiplier \eqref{eq3.3}. 



The final remark we want to make about Theorem \ref{3.T.1} is that the condition on the initial datum is not so easy to verify.  In fact to check whether or not the condition is satisfied is quite a delicate question in Fourier Analysis. In the case of the free Schr\"odinger operator this question is nothing but to study $\,L^p\,$ bounds for the dual operator to the free propagator $\,e^{it\Delta}u_0\,$. In terms of this latter operator the question is if the inequalities                                      
\begin{equation}                                                                                                                            
\label{eq3.5}                                                                                                                               
\left\|e^{it\Delta}u_0\right\|_{L^p(\mathbb{R}^{d+1})}\le c\left\|\widehat{u}_0\right\|_{L^q(\mathbb{R}^d)}                                 
\end{equation}                                                                                                                              
with                                                                                                                                        
\begin{itemize}                                                                                                                             
 \item $\D\frac{d+2}p+\D\frac dq=d$                                                                                                         
 \item $p>d$,                                                                                                                               
\end{itemize}                                                                                                                               
hold true. The answer is positive if $\,d=1\,$ and it was established by C. Fefferman and E. Stein (see \cite{Fe}). In higher dimensions $\,d\ge 2\,$ this was proved by E. Stein and P.A. Tomas for $\,q\le 2\,$ (see \cite{St}). Recently, the case $d=2$ was extended by T. Tao in \cite{T} to $\,q\le\D\frac 52\,$ . As a consequence if $\,\mu\,$ is small enough and
 \begin{equation}
\label{eq3.6}
\left\|\widehat{u}_0\right\|_{L^{q,\infty}}\le \mu
\end{equation}
for $\,q=\D\frac{d\alpha}{d\alpha-2}\,$ with $\,\D\frac 83<\alpha<8\,$ if $\,d=1\,$ and $\,\frac 53<\alpha<3\,$ if $\,d=2\,$ we can apply Theorem \ref{3.T.1}. Moreover in \cite{CVV} it is proved that  if $\,u_0\,$ verifies \eqref{eq3.6} then 
$$\sup_t\left\|\widehat{u}(\cdot,t)\right\|_{L^{q,\infty}}\le C\mu$$
for a fixed constant $\,C\,$ and with $\,q\,$ as above.

This was the starting point of the work by A. Gr\"unrock \cite{Gr1}\cite{Gr2}. In \cite{Gr2} he proves the following theorem
\begin{theorem}\label{3.P.2} Let $\,1<r\le 2\,$ and $\,\widehat{u}_0\in L^{r^{\prime}}(\R)\,$. Then the I.V.P.
\begin{equation}
\label{eq3.7}
iu_t+u_{ss}=\lambda |u|^2u\qquad u(0)=u_0\quad \lambda=\pm 1,
\end{equation}
is locally well posed and
\begin{equation}
\label{eq3.8}
\widehat{u}(\cdot,t)\in L^{r^{\prime}}(\R)
\end{equation}
for all $t$.
\end{theorem}

So in particular if
$$\left|\widehat{u}_0(\xi)\right|\le\D\frac{1}{(1+|\xi|)^{\gamma}},\qquad \gamma>0,$$
the conclusion holds. As it shall be explained later, the selfsimilar solutions of BF will lead to the case $\,\gamma=0$.

Finally we would like to mention another work, by A. Vargas and L. Vega \cite{VV}, related to \eqref{eq3.7} and which precedes to \cite{Gr2}. First, the following local result is proved in \cite{VV}.

\begin{proposition}\label{3.P.3} Assume $\,u_0\in \mathcal{S}^{\prime}(\R)\,$ such that there exists $\tilde{T}>0$ with
\begin{equation}
\label{eq3.9}
\left\|e^{it\partial^2_x} u_0\right\|_{L^3([-\widetilde{T},\widetilde{T}],L^6(\mathbb{R}))}<+\infty.
\end{equation}
Then there exists $\,T\le \widetilde{T}\,$ and a unique solution $\,u\,$ of \eqref{eq3.7} such that
$$\|u\|_{L^3([-T,T],L^6(\mathbb{R}))}<+\infty.$$
\end{proposition}

We say that $\,u_0\in \mathcal{Y}_2\,$ if \eqref{eq3.9} is satisfied uniformly in $\tilde{T}$. As a global result, the following theorem is proved in \cite{VV}. 

\begin{theorem}\label{3.T.2} Let $\,u_0\in \mathcal{S}^{\prime}(\R)\,$ that can be decomposed as $\,u_0=v_0+w_0\,$ with
\begin{equation}
\label{eq3.10}
\|v_0\|_{L^2(\mathbb{R})}\le M^{r}\quad,\quad \|w_0\|_{\mathcal{Y}_2}\le M^{-1}\quad\mbox{for some }\,r<1\,\mbox{ and  }\,M>1.
\end{equation}
Then there exists a unique solution $u$ of \eqref{eq3.7} for
$$0\le t\le T= T(M),$$
that satisfies
$$\|u\|_{L^3([0,T],L^6(\mathbb{R}))}\le c.$$
Moreover
$$\lim_{M\rightarrow \infty}T(M)=\infty.$$
\end{theorem}

As a consequence if for all $\,M>0\,$, $u_0$ can be decomposed as $\,u_0=v_0^M+w_0^M\,$ satisfying \eqref{eq3.10}, then there is a unique global solution of \eqref{eq3.7}.

Again the difficulty relies in giving good conditions on $\,u_0\,$ such that \eqref{eq3.9} is satisfied. In terms of the pointwise decay of the Fourier transform it is worse than the result by Gr\" unrock, because a decay $r>1/6$ is necessary. However in terms of the $L^p$ local integrability of $\hat u_0$ the results are stronger. In fact it is sufficient to assume the following. Assume $\,\widehat{u}_0\in X_{p,\gamma}\,$ which is given as the functions $\,f\,$ such that
  $$\|f\|_{p,\gamma}=\sum_{j\ge 0}
  \left(\sum_{\mathcal{T}\in \mathcal{D}_j} 2^{-p\gamma j}\|f_{|_\mathcal{T}}\|^p_{L^2}\right)^{1/p}$$
  where
  $\,\mathcal{D}_j=\{[k2^j,(k+1)2^j)\,:\,k\in\Z\}\,$, the family of dyadic intervals of length $\,2^j\,$ in $\,\R\,$. Then \eqref{eq3.9} holds $\,p=3\,,\,\gamma=\D\frac 16$. The above norms were firstly introduced by J. Bourgain and naturally appear  in the study of the free solutions of Schr\"odinger equation - see also \cite{PlV} equations, (3.1) and (3.2).
   
 Notice that the advantage of Proposition \ref{3.P.3} with respect to Theorem \ref{3.P.2}  is that the solution obtained in the former case is $\,L^1_{loc}(\mathbb{R}^{d+1})$. This property is fundamental in order to integrate the Frenet equations to obtain a solution of the BF from cubic NLS. In any case it would be interesting to know if the methods by A. Vargas and L. Vega can be extended to the same range of decay as the ones obtained by A. Gr\"unrock.

 \subsection{The $\,\delta$--function as initial datum for NLS}
Assume next that $\,u_0=a\delta\,$, for $a\in\mathbb{C}$. Then for all $N$
$$u_{0,N}(s):=e^{iNs}u_0=u_0.$$
If uniqueness holds then we conclude from the galilean transform \eqref{eq1.2} that
$$u(s,t)=u_N(s,t)=e^{-itN^2+iNs}u(s-2Nt,t)\quad\mbox{ for all }\,N.$$
After differentiation with respect to $N$ in the above identity we get
$$0=i(s-2Nt)\,u(s-2Nt,t)-2t\,\partial_su(s-2Nt,t).$$
Hence
$$u(s,t)=g(t)e^{i\frac{s^2}{4t}}.$$
But $\,u\,$ has to solve NLS. This gives
\begin{equation}
\label{eq1.3}
g(t)=\D\frac{a}{\sqrt t}\exp \left(\pm i|a|^2\ln t+ia_1\right),\quad a,a_1\in\C,
\end{equation}
which has no limit as $\,t\uparrow 0^+\,$. Therefore (\ref{eq1.1}) is ill posed with $\,u_0=a\delta$ (see the article of Kenig, Ponce and Vega \cite{KPV2}).

However notice that if we change NLS from (\ref{eq1.1}) into
$$iu_t+u_{ss}+\D\frac 12\left(|u|^2+A(t)\right)u=0,$$
we see that taking $\,A(t)=-\D\frac{|a|^2}{t}\,$, then $\,u_{a}(s,t)=\D a\frac{e^{i\frac{s^2}{4 t}}}{\sqrt t}\,$ is a solution of
$$\left\{\begin{array}{rcl}
iu_t+u_{ss}+\D\frac 12\left(|u|^2-\D\frac{|a|^2}{t}\right)u&=&0,\\
u(s,0)&=&a\delta.
\end{array}\right.$$

This suggests the possibility of finding one solution of the binormal flow. The corresponding curve should be determined by the property that curvature and torsion are
\begin{equation}
\label{eq1.4}
c(s,t)=\D\frac{a}{\sqrt t}\qquad,\qquad
\tau(s,t)=\D\frac{s}{2t},
\end{equation}
for some $a\in\mathbb{R}^+$.

\section{Selfsimilar solutions for the binormal flow}
In this section we review some of the results proved in \cite{GuRV} and  
\cite{GV} concerning the self-similar solutions of the binormal flow.  We 
also include in Propositions  \ref{L0} and \ref{L00} some  results for curves in  
$\R^3$
that as far as we know are new. These results are an extension to  general 
curves of those proved in \cite{GuRV} for the particular case  in which the 
curve is precisely a self-similar solution of the  binormal 
flow.
\subsection{A V-shape as initial datum for the binormal flow}

The BF solution corresponding to (\ref{eq1.4}) should have as initial condition
$$\chi(s,0)=\begin{cases}
      sA^+& \text{if }s\geq0 \\
      -sA^-& \text{if }s<0,
\end{cases}$$
with $\,|A^+|=|A^-|=1$.

Frenet equations \eqref{Frenet} allow us to compute $\,\chi(s,t)\,$ from (\ref{eq1.4}) up to a rigid motion, obtaining a real analytic curve for $\,t>0\,$ which should have a singularity at $\,t=0\,$. Notice that the binormal flow is reversible in time: if $\chi(s,t)$ is a solution so is $\tilde\chi(s,t)=-\chi(s,-t)$. Recall that the energy in this case is infinity, although it is locally finite as for the helix.

We gather the results proven on selfsimilar solution of the binormal flow in the following therorem.
\begin{theorem} The self-similar solutions are parametrized by a positive number $a$, and satisfy 
$$\left|\chi_a(s,t)-sA^+_a  \mathbb{I}_{[0,\infty)}(s)-sA^-_a  \mathbb{I}_{(-\infty,0]}(s)\right|\leq {2a}{\sqrt{t}},$$
for $t\geq 0$, where  $A^\pm_a$ are two distinct non opposite unit vectors. 
In particular, 
\begin{equation}
\chi_a(s,0)=\begin{cases}
     A^+_a s & s>0 \\
    A^-_a s & s<0.
\end{cases}
\end{equation}
Moreover, there exists $\tilde{a},\tilde{\tilde{a}}\,$ such that if $\,a<\tilde{a}\,$, $\,\chi_{a}(s,1)\,$ has no self-intersection and if $\,a>\tilde{\tilde{a}}\,$ it has finitely many.
Finally, the relation between $a$ and $A^\pm_a$ is
$$\sin\frac\theta 2=e^{-\frac{\pi a^2}{2}},$$
where $\theta$ is the angle between the vectors $A^+_a$ and $-A^-_a$.
\end{theorem}
The behaviour of $\chi_a(s,t)$ for $t>0$ was exhibited by the physicists M. Lakshmanan and M. Daniel in \cite{LD} and a numeric 
study on it was done 
by T.F. Buttke in \cite{Bu1}. The rigorous description of these solutions up to $t=0$ was done by S. Guttierez, J. Rivas and the second 
author in \cite{GuRV}, and will be presented in the two next sections. Let us mention that the link between the parameter $a$ and
the angle 
between the vectors $A^+_a$ and $-A^-_a$ is an original material, and we dedicate to its proof the third subsection  \S 3.3.

Before trying to compute $\,\chi\,$ is worthy to stress the point that galilean transformations act on the set of solutions of the BF. In fact writting the flow in intrinsic coordinates we get -see \cite{DaR}, \cite{Be}
\begin{equation}
\label{eq1.41}
\begin{array}{l}
c_t=-(c\tau)_s-c_s\tau\\
\tau_t=\left(\D\frac{c_{ss}-c\tau^2}{c}\right)_s+c_sc.
\end{array}
\end{equation}

It is straighforward to observe that for $\,N\in\R$
$$\left(c_N,\tau_N\right)=
\left(c(s-2Nt,t),\tau(s-2Nt,t)+N\right)$$
is a solution of (\ref{eq1.41}) if so is $(c,\tau)$. Also 
$$\left(c_N,\tau_N\right)=(c,\tau)\quad\forall N\qquad\mbox{iff}\qquad
(c,\tau)=\left(\D\frac{a}{\sqrt t},\D\frac{s}{2t}\right)\quad t>0.$$

\subsection{The computation of the selfsimilar solutions and some remarks}

The simplest way to construct $\chi(s,t)$ is to observe that from (\ref{eq1.4}) we can expect it to have a selfsimilar behaviour also with respect to scaling. Hence we settle the ansatz:
$$\chi(s,t)=\sqrt t \,G\left(\frac{s}{\sqrt t}\right)\quad t>0,$$
which gives from the BF equation ($\,\chi_t=\chi_s\land\chi_{ss}\,$)
\begin{equation}
\label{eq1.5}
\D\frac 12G-\frac s2G^{\prime}=G^{\prime}\land G^{\prime\prime}.
\end{equation}

Since $\,G^{\prime}=T\,$, a new differentiation gives using Frenet equations 
$$\begin{array}{rcl}
-\D\frac s2\,T^{\prime}=T\land T^{\prime\prime}&=&
\big(T\land (c_s\,n+c\,n_s)\big)\\
&=&c_s\,b+c\,\big(T\land (-c\,T+\tau \,b)\big).
\end{array}$$
Hence (see \cite {LD} and \cite{Bu1}),
\begin{equation}
\label{eq1.6}
\D\frac s2\,c\,n=-c_s\,b+c\,\tau \,n
\end{equation}
and therefore
\begin{equation}\quad c_s=0\quad \tau=\frac s2.\end{equation}
In particular, there exists $a\in\mathbb{R}^+$ such that $c=a$.
Now the Frenet equations write
\begin{equation}
\label{eq1.7}
\begin{array}{rclll}
T'&=&&a\,n\,,\\
n'&=&-a\,T&&+\,\D\frac s2\,b\,,\\
b'&=&&-\,\D\frac s2\,n\,,
\end{array}
\end{equation}
and we get that $\,G^{\prime}=T\,$ is a real analytic function. Moreover
\begin{equation}
\label{eq1.8}
T(s,t)=T\left(\frac{s}{\sqrt t}\right).
\end{equation}

Recall that BF is invariant under rotations so we can assume without loss of generality that
$$\left(\begin{array}{c}
T\\n\\b\end{array}\right)(0,0)=\mathbb{I}_{3\land 3}\quad\mbox{ and from (\ref{eq1.8}) }\quad
\left(\begin{array}{c}
T\\n\\b\end{array}\right)(0,t)=\mathbb{I}_{3\land 3}.$$

We have to determine $\,\chi(0,t)\,$. But from
$$\chi_t=cb$$
we conclude
$$\chi_t(0,t)=\D\frac{a}{\sqrt t}\,(0,0,1).$$

As a conclusion $\,\chi(0,t)=2a\sqrt t\,(0,0,1)\,$. We could have got it also from (\ref{eq1.5}).

Therefore
$$\chi(s,t)=\sqrt t \,G\left(\frac{s}{\sqrt t}\right)$$
and $\,G\,$ is determined up to a rotation by the solution of (\ref{eq1.7}) with
$$\begin{array}{l}
G(0)=2a\,(0,0,1)\\
T(0)=(1,0,0)\\
n(0)=(0,1,0).
\end{array}$$

It is easy to do the numerical analysis of $\,G\,$. We start with \eqref{eq1.5}:
$$
\D\frac 12G-\frac s2G^{\prime}=G^{\prime}\land G^{\prime\prime}$$
Recall the identity:
\begin{equation}
\label{eq1.9}
F\land(G\land H)=(F\cdot H)\cdot G-(F\cdot G)\cdot H.
\end{equation}
Then
$$G^{\prime}\land\left(G^{\prime}\land G^{\prime\prime}\right)=
\left(G^{\prime}\cdot G^{\prime\prime}\right)G^{\prime}-
\left(G^{\prime}\cdot G^{\prime}\right)G^{\prime\prime}$$
But on the other hand from
$$-\D\frac s2\,T^{\prime}=T\land T^{\prime\prime}$$
we get
$$-\D\frac s2\, T^{\prime}\cdot T=0=-\D\frac s2\,G^{\prime}\cdot G^{\prime\prime}$$
and $\,|T|^2=1$.

Hence from \eqref{eq1.5}, \eqref{eq1.9}
\begin{equation}
\label{eq1.10}
\left\{\begin{array}{rcl}
\D\frac 12G\land G^{\prime}&=&G^{\prime\prime}\\
G^{\prime}(0)&=&(1,0,0)\qquad\qquad G(0)\cdot G^{\prime}(0)=0,
\end{array}\right.
\end{equation}
which can be easily computed numerically. We refer the reader to the web page http://ehu.es/luisvega/vortex/files/esquina.mov to see 
the corresponding flow of curves. In these pictures is clearly seen that the profile $\,G\,$ tends to two straight lines at infinity, and as a consequence that
\begin{equation}
\label{eq1.11}
\chi(s,0)=\begin{cases}
     A^+s &s>0 \\
    A^- s & s<0.
\end{cases}
\end{equation}

In order to prove it we will have to work a little bit. We follow the arguments of S. Gutierrez, J. Rivas and L. Vega \cite{GuRV}. Recall that from \eqref{eq1.5}:
$$\D\frac 12G-\frac s2\,G^{\prime}=T\land T^{\prime}=T\land a\,n=a\,b.$$

Hence, for $\,s>0$
$$\left(\D\frac {G(s)}{s}\right)^{\prime}=\D\frac{s\,G^{\prime}-G}{s^2}=-\D\frac{2a}{s^2}\,b.$$
It follows that
\begin{equation}
\label{eq1.12}\D\frac{G(s)}{s}=G(1)+\int_1^s-\D\frac{2a}{s^{\prime2}}\,b(s')\,ds',
\end{equation}
and therefore, since $b$ is unitary, we conclude that there exists 
\begin{equation}
\label{eq1.13}\,\D\lim_{s\to\infty}\D\frac{G(s)}{s}=A^+.
\end{equation}
Analagously we obtain that there exists $\,A^-\,$ such that $\D\lim_{s\to-\infty}\D\frac{G(s)}{s}=A^-.$

Moreover,
$$A^+-\D\frac{G(s)}s=\left.\D\frac{G(s)}s\right|_s^{\infty} =-2a\int_s^{\infty}\D\frac {b(s')}{{s^{\prime}}^2}\,ds',$$
and 
   \begin{eqnarray}
\label{eq1.14} G(s) & = & sA^++2as\D\int_s^{\infty}\D\frac {b(s^{\prime})}{{s^{\prime}}^2}ds^{\prime}\qquad s>0 \\
\label{eq1.15} G(s) & = & sA^--2as\D\int_{-\infty}^{s}\D\frac {b(s^{\prime})}{{s^{\prime}}^2}ds^{\prime}\qquad s<0.
\end{eqnarray}
We have the following proposition.
\begin{proposition}\label{P1} For all $s$
$$\left|\chi(s,t)-sA^+\mathbb{I}_{[0,\infty)}(s)-sA^-\mathbb{I}_{(-\infty,0]}(s)\right|
\le 2a\,\sqrt t\qquad t>0,$$
where $\mathbb{I} $ denotes  the characteristic function.
\end{proposition}
\underline{Proof}
If $s=0$ it is trivial because  $\, \chi(0,t)=2a\sqrt t$. We will assume that $s>0,$ the other case been proved in a similar way. Then
$$
\chi(s,t)=\sqrt tG\left(\D\frac s{\sqrt t}\right)=sA^++2as\D\int_{s/\sqrt t}^{\infty} \,\D\frac{b(s')}{{s^{\prime}}^2}\,ds^{\prime},
$$
and then
$$\left|\chi(s,t)-sA^+\right|\le 2as
\left|\D\int_{s/\sqrt t}^{\infty}\,\D\frac{b\left(s^{\prime}\right)}{{s^{\prime}}^2}\,ds'\right|
.$$
The proof follows using that $b$ is unitary.

\begin{remark}
It can be also proved that (see  \cite{GuRV})

$$\D\lim_{t\to 0} \,T(s,t)=\begin{cases}
      A^+& s>0, \\
     A^- & s<0.
\end{cases}$$
\end{remark}

Our next purpose is to prove that $\,A^+\neq A^-\,$. Recall that we are dealing with a one parameter family of solutions with $a$ as the free parameter. We will make explicit the dependence on $a$
in the following proposition in order to clarify the exposition.

\begin{proposition}\label{P2}
\begin{itemize}
  \item [(i)] $A^+_a\neq -A^-_a\quad \forall a.$
  \item [(ii)] $A^+_a=A^-_a\quad\Longleftrightarrow\quad a=0$.
  \item [(iii)] $\exists\, \,\tilde{a},\tilde{\tilde{a}}\,$ such that if $\,a<\tilde{a}\,$, $\,G_{a}(s)\,$ has no selfintersection and if $\,a>\tilde{\tilde{a}}\,$ it has finitely many.
\end{itemize}
\end{proposition}
\underline{Proof}  Recall that $\,G_{a}-s\,G_{a}^{\prime}=2a\,b$. Also $\,G_{a}^{\prime\prime}=T^{\prime}=a\,n$. Hence
$$\begin{array}{rcl}
G_{a}^{\prime\prime\prime}=an^{\prime}&&=a\left(-a\,T+\D\frac s2\, b\right)\\\\
&&=-a^2\,T+\D\frac s4\left(G_{a}-s\,G_{a}^{\prime}\right).
\end{array}$$
Then
\begin{eqnarray}
\label{eq1.16}
G_{a}^{\prime\prime\prime}+\left(a^2+\D\frac{s^2}4\right)G_{a}^{\prime}-\D\frac s4\,G_{a}&=&0\\\nonumber\\
\left\{\begin{array}{rcl}
G_{a}(0)&=&2a\,(0,0,1)\\
G_{a}^{\prime}(0)&=&(1,0,0)\\
G_{a}^{\prime\prime}(0)&=&a\,(0,1,0)
\end{array}\right.
\label{eq1.17}
\end{eqnarray}
Also,
\begin{itemize}
  \item $G_{a}(s)=\big(x_{a}(s),y_{a}(s),z_{a}(s)\big)$  with $\left(\begin{array}{c}
  x_{a}(0)\\x_{a}^{\prime}(0)\\x_{a}^{\prime\prime}(0)
  \end{array}\right)=\left(\begin{array}{c}
  0\\1\\0
  \end{array}\right)$.
  \item In \eqref {eq1.16} if $G_{a}(s)$ is a solution then $G_{a}(-s)$ is a solution.
  \end{itemize}
We conclude that if $a\neq 0$ then
\begin{eqnarray}
\begin{array}{rcl}
x_{a}(s)&=&-x_{a}(-s),\\
y_{a}(s)&=&y_{a}(-s),\\
z_{a}(s)&=&z_{a}(-s),
\end{array}
\label{eq1.18}
\end{eqnarray}
and therefore $x_{a},\,y_{a},\,z_{a}$ are three linearly independent solutions of the scalar ODE
\begin{equation}
\label{eq1.19}
w^{\prime\prime\prime}+\left(a^2+\D\frac{s^2}4\right)w^{\prime}-\D\frac s4\, w=0
\end{equation}
which are tempered distributions $(G_{a},G_{a}^{\prime},G_{a}^{\prime\prime},G_{a}^{\prime\prime\prime}\mbox{ grow polynomialy})$. Applying the Fourier transform we get:
\begin{equation}
\label{eq1.20}
\xi\,\widehat{\omega}^{\prime\prime}+3\,\widehat{\omega}^{\prime}+4\xi^2\,\widehat{\omega}-
4\xi a^2\,\widehat{\omega}=0.
\end{equation}
This is a singular regular O.D.E. with $r=0$ and $r=-2$ the solutions of the indicial equation
$$r(r-1)+3r=0.$$

Therefore the regular solutions have to be even, and the odd solutions have to behave as $\omega_0\,\delta^{\prime}$ at the origin. But we know that $x_{a}$ is odd and non trivial. Therefore $w_{a}\neq 0$ which amounts to say that $A^+_{a,1}\neq 0$, where we use the vectorial notation $A^\pm_a=(A^\pm_{a,1},A^\pm_{a,2},A^\pm_{a,3})$. From \eqref{eq1.18} we get
$$\begin{array}{rcl}
A_{a,1}^+&=&A_{a,1}^-,\\
A_{a,2}^+&=&-A_{a,2}^-,\\
A_{a,3}^+&=&-A_{a,3}^-.
\end{array}$$
Hence $A^+_a+A^-_a\neq 0$.

Assume now that $A^+_a=A^-_a$. Then $A_{a,2}^+=A_{a,2}^-=A_{a,3}^-=A_{a,3}^+=0$. But then the even functions $\widehat y_{a}(\xi),\,\widehat z_{a}(\xi)$ should be regular at the origin and linearly dependent. Hence $a=0$ and $A^+_a=(1,0,0)=A^-_a$.

Only (iii) remains to be proved. Recall that $G_{a}=s\,T+2a\,b$. Hence $|G_{a}|^2=s^2+a^2$. Therefore if $G_{a}(s_1)=G_{a}(s_2)$ for different $s_1$ and $s_2$, then $\,s_1^2=s_2^2$ and $s_1=-s_2$. But from \eqref{eq1.18} this implies that $x(s_1)=0=x(s_2)$. Recall that
$$\begin{array}{rcl}
x^{\prime\prime\prime}_{a}+\left(a^2+\D\frac{s^2}4\right)x^{\prime}_{a}-\D\frac s4\, x_{a}&=&0,\\
x_{a}(0)&=&0,\\
x_{a}^{\prime}(0)&=&1,\\
x_{a}^{\prime\prime}(0)&=&0.
\end{array}$$

Hence $x_{a}(s)=s+O\left(a^2\,s^3\right)$ and on the other hand the map $a\longmapsto \lim\limits_{s\to\infty}\D\frac{x_a(s)}s=A_{1,a}^+$ is continuous from the proof of Proposition \ref{P1}. Hence if $a$ is small enough there is no selfintersection. On the other hand it is easy to see that 
$x^{\prime}_{a}(s)=\cos(as)+O\left(\D\frac{s^3}{a}\right)$ and therefore
$$x_{a}(s)=\D\frac 1{a}\sin (as)+O\left(\D\frac{s^4}{a}\right)\quad 0<s<1.$$
Hence for $a$ large there are selfintersections.

Notice that if $\chi(s,t)$ is a solution $-\chi(s,-t)$ is also a solution. Hence
$$\left\{\begin{array}{l}
\chi_t=\chi_s\,\land \chi_{ss}\\
\chi(s,-1)=-G_{a}(s).
\end{array}\right.$$
is a solution that develops a singularity at time $t=0$ although $G_{a}$ is real analytic.

A few remarks are in order. 
\begin{remark}
The kinetic energy $\left(\int c^2(s,1)\,ds\right)$ is infinity but locally finite as in the case of the helix.
\end{remark}

\begin{remark} The solution $\chi_{a}(s,t)=\sqrt t \,G_{a}\left(\D\frac s{\sqrt t}\right)$ has as torsion $\D\frac s{2t}$ which is selfsimilar of the type
  $$\tau_{\lambda}(s,t)=\lambda \,\tau(\lambda s,\lambda^2 t).$$
  This scaling leaves invariant the following conserved quantity
  $$\int \tau (s,0)\,ds=\int\tau (s,t)\,ds.$$
  Notice that in our situation the above invariant has to be understood in the sense
  $$\lim\limits_{R\to\infty}\int_{-R}^R \tau (s,t)\,ds=0.$$
\end{remark}

\begin{remark} It is interesting to analyze the type of singularity we get from the point of view of Euler equations and vortex filaments. In particular let us recall the condition obtained by P. Constantin, C. Fefferman and A. Majda \cite {CFM} which is  measured in terms of $\|\nabla_x\xi \|_{L^{\infty}}$, where $\xi$ stands for the unit direction of the vorticity:
  \begin{equation}
\label{eq1.19bis}
\xi=\D\frac{\omega}{|\omega|}.
\end{equation}
Their condition is that
\begin{equation}
\label{eq1.20bis}
\int_0^T\|\nabla_x\xi\|^2_{L^{\infty}}dt=+\infty
\end{equation}
if there is a blow up at time $T$. In our setting $\xi$ plays the role of $T$ and therefore $|\nabla_x\xi|$ is the curvature. Hence \eqref{eq1.20bis} becomes
\begin{equation}
\label{eq1.21}
\int_{-1}^0\D\frac{a^2}{|t|}\,dt=\infty.
\end{equation}
\end{remark}

\begin{remark} Recall that in the case of the free Schr\"odinger equation we have a similar situation. Taking $$u(s,t)=\D\frac {e^{i\frac{s^2}{4t}}}{\sqrt{|t|}}\qquad -1<t<0,$$
  we have a solution of
  $$\left\{\begin{array}{rcl}
  iu_t+u_{ss}&=&0,\\
  u(s,-1)&=&e^{i\frac{s^2}{4}},
  \end{array}\right.$$
  such that $u(s,0)=2\sqrt{\pi}e^{i\frac{\pi}{4}}\,\delta$.
\end{remark}

\begin{remark}
Equation \eqref{eq6} of the tangent vector suggests many possible generalizations by considering other targets besides the sphere. 
Let us then introduce the notation
$$u\land_- v=\left(\begin{array}{ccr}
1&0&0\\
0&1&0\\
0&0&-1
\end{array}\right)u\land v.$$
Therefore instead of \eqref{eq6} we consider
\begin{equation}\label{mapsign}
T_t=T\land_- T_{ss},
\end{equation}
with $T$ a map from $\R^2$ onto the hyperbolic plane $\mathbb{H}^2$. Analogously, since  $T=\chi_s$, we can obtain the equation 
\begin{equation}\label{binormalsign}
\chi_t=\chi_s\land_-\chi_{ss}.
\end{equation}
Similar calculations and results as in \cite{GuRV} were done by de la Hoz \cite{Pa} for selfsimilar solutions of \eqref{binormalsign} in the hyperbolic setting. The extra-difficulty is that there is a-priori no control on the size of the euclidean length of the generalized binormal vectors. \end{remark}

\begin{remark}                                                                                                                              
Let us finally remark that numerical computations on the selfsimilar curves of BF were done by F. de la Hoz,  C. Garc\'ia-Cervera and L. Vega in \cite{HGCV}. In that paper the authors point out the similarity of these curves with those that appear in the experiments of the flow around a delta wing.                                                                                                                   
\end{remark}

\subsection{Proof of the identity $A_{a,1}=e^{-\frac{\pi a^2}{2}}$}

Let us try to answer a more delicate question. We have been able to construct non trivial solutions of BF such that  the initial condition has a V-shape determined by two unit vectors $A^+_a$ and $A^-_a$. The question we propose is if given any pair of vectors $(A^+,A^-)$ there is $a$ such that 
$(A^+,A^-)=(A^+_a, A^-_a)$.
 We have seen in Proposition \ref{P2} (ii), that the particular case $A^+=-A^-$ cannot be obtained using the family of curves $\{G_{a}\}$, ${a\in\R}$. In this section we shall prove that except for this case the answer to our question is positive.
 
  The proof requires to integrate the Frenet system of equations in a different way.
  We start defining the complex vector
  $${N}={n}+i\,{b}\,e^{i\int_0^s\tau(s^{\prime})ds^{\prime}}$$
  and
  $$u=c\,e^{i\int_0^s\tau(s^{\prime})ds^{\prime}}.$$
  Then we rewrite the Frenet system as:
  $$\begin{array}{rcl}
  {N}_s&=&-u\,{T}\\
  {T}_s&=&\D\frac 12\left(\overline{u}\, {N}+u\,{\overline{N}}\right).
  \end{array}$$
  
  We denote $(v_1,v_2,v_3)$ the components of a $\mathbb{R}^3$ vector $v$. Then
  $$\left|n_j+ib_j\right|^2+T_j^2=n_j^2+b_j^2+T_j^2=1.$$
  Next we introduce
  $$\varphi_j=\D\frac{N_j}{1-T_j}$$
so that $\varphi_j$  satisfies the Riccati equation
  $$\varphi_j^{\prime}-\D\frac 12\,\overline{u}\,\varphi_j^2-\D\frac 12\,u=0.$$
  Now we consider $\eta_j$ such that
  $$\begin{array}{c}
\varphi_j(s)=\D\frac{u}c\,\eta_j(s)=e^{i\int_0^s \tau}\,\eta_j,\\\\
\varphi_j^{\prime}=\D\frac{u}c\eta_j^{\prime}+i\tau\D\frac{u}c\eta_j=
\D\frac 12u\eta_j^2+\D\frac 12u=\D\frac 12u\left(\eta_j^2+1\right).
\end{array}$$
Hence
$$\eta_j^{\prime}+i\tau\eta_j=\D\frac 12c\left(1+\eta_j^{2}\right).$$

The next step is to define $\theta_j$ as $\theta_j=\exp[-c/2\int_0^s\eta_j]$, i.e.:
$$\eta_j=-\D\frac 2c\,\frac{\theta_j^{\prime}}{\theta_j}$$
Then
$$\eta_j^{\prime}=2\D\frac{c^{\prime}}{c^2}\,\frac{\theta_j^{\prime}}{\theta_j}
-\D\frac 2c\,\frac{\theta_j^{\prime\prime}}{\theta_j}+
{\D\frac 2c\,\frac{{\theta_j^{\prime}}^2}{\theta_j^2}}$$
Hence
$$
2\D\frac{c^{\prime}}{c^2}\,\frac{\theta_j^{\prime}}{\theta_j}-
\D\frac 2c\,\frac{\theta_j^{\prime\prime}}{\theta_j}=-i\tau\D\frac 2c\, \frac{\theta_j^{\prime}}{\theta_j}+
\D\frac 12\, c.$$
We have arrived at the remarkable fact that $\theta_j$ solves
\begin{equation}
\theta_j^{\prime\prime}+\left(-\D\frac{c^{\prime}}c+i\tau\right)\theta_j^\prime
+\D\frac{c^2}4\theta_j=0\label{eq1.23}.
\end{equation}
Multiplying by $\,\D\frac{\theta^{\prime}_j}{c^2}\,$ we get that
  \begin{equation}
\label{eq1.24}
E(s):=\left(\D\frac{\theta^{\prime}_j}c\right)^2+\D\frac{\theta^{2}_j}4=\mbox{constant}=E(0).
\end{equation}
Also notice that
$$|N_j|^2+|T_j|^2=1$$
is equivalent to $(n_j-ib_j)(n_j+ib_j)=(1-T_j)(1+T_j)$. Hence
$$|\varphi_j|^2=\varphi_j \overline{\varphi_j}=\D\frac{n_j+ib_j}{1-T_j}\cdot 
\D\frac{n_j-ib_j}{1-T_j}=\D\frac{1+T_j}{1-T_j},$$ and
$$T_j=\D\frac{|\varphi_j|^2-1}{|\varphi_j|^2+1}.$$
But $\varphi_j=e^{i\int _0^s\tau}\eta_j$ and $\eta_j=-\D\frac{2\theta_j^{\prime}}{c\theta_j}$. Then
$$T_j=\D\frac{|\eta_j|^2-1}{|\eta_j|^2+1}=\D\frac{\frac 4c|\theta_j^{\prime}|^2-|\theta_j|^2}{\frac 4c|\theta_j^{\prime}|^2+|\theta_j|^2}=1-\D\frac{|\theta_j|^2}{2E(0)}.$$

Therefore we have proved the following proposition valid for any curve.

\begin{proposition}\label{L0} Any of the components $T_j$ of the tangent vector $T$ can be written as
$$T_j=1-\D\frac{|\theta_j|^2}{2E(0)},$$
with $\theta_j$ solution of (\ref{eq1.23}) and $E(0)$ as in (\ref{eq1.24}).
\end{proposition}

Assume for simplicity
$$\left(\begin{array}{c}
T\\n\\b\end{array}\right)(0)=\mathbb{I}_{3\land 3}.
$$
Then
$$T_1(0)=1\quad\Longrightarrow\quad \theta_1(0)=0.$$
Taking $\,\theta_1^{\prime}(0)=\frac{c(0)}{\sqrt 2}\,$ we get
  \begin{equation}
\label{eq1.25}
T_1=1-|\theta_1|^2\quad;\quad E(0)=\frac 12.
\end{equation}
Similarly 
$$T_2(0)=0\quad\Longrightarrow\quad\eta_2(0)=1\quad,\quad
\theta_2(0)=1\quad,\quad \theta_2^{\prime}(0)=-\D\frac{c(0)}2.$$
Then
 \begin{equation}
\label{eq1.26}
T_2=1-|\theta_2|^2\quad;\quad E(0)=\frac 12.
\end{equation}
Finally
$$T_3(0)=0\quad,\quad\eta_3(0)=i\quad,\quad
\theta_3(0)=i\quad,\quad \theta_2^{\prime}(0)=\D\frac{c(0)}2$$
gives
 \begin{equation}
\label{eq1.27}
T_3=1-|\theta_3|^2\quad;\quad E(0)=\frac 12.
\end{equation}
As for $n$ and $b$ we get, dropping the $j$--index,
\begin{itemize}
  \item [$\bullet$] $cn=T_s=\left(1-|\theta|^2\right)_s=-2\mbox{Re}\,\theta\overline{\theta}^{\prime}$
  \item [$\bullet$] $T_{ss}=-2|\theta^{\prime}|^2-2\mbox{Re}\,\theta\overline{\theta}^{\prime\prime}=-2|\theta^{\prime}|^2+\D\frac{c^2}2|\theta|^2-2\D\frac{c^{\prime}}c\mbox{Re}\,\theta\overline{\theta}^{\prime}+2\tau\mbox{Im}\,\theta\overline{\theta}^{\prime}$
   \end{itemize}
  From $\,T_{ss}=(cn)_s\,$ we get
  $$\D\frac c2b=\mbox{Im}\,\theta\overline{\theta}^{\prime}.$$
  Hence we have proved the following proposition valid for any curve.
  \begin{proposition}\label{L00} Any of the components $n_j$ and $b_j$  of the normal and binormal vectors $n$ and $b$ can be written as
   \begin{equation}
\label{eq1.28}
c(n_j-ib_j)=-2\theta_j\overline{\theta}^{\prime}_j\quad,\quad j=1,2,3,
\end{equation}
with $\theta_j$ solution of (\ref{eq1.23}) verifying (\ref{eq1.25})-(\ref{eq1.27}).

\end{proposition}

In our particular case equation \eqref{eq1.23} reduces to
\begin{equation}
\label{eq1.29}
\theta^{\prime\prime}+i\D\frac s2\,\theta^{\prime}+\D\frac{a^2}4\,\theta=0.
\end{equation}
Our final purpose is to find a relation between $\gamma$, the angle determined by $A^+_a$ and $-A^-_a$, and the parameter $a$.

It is easy to find two linearly independent solutions by just computing the Fourier transform in \eqref{eq1.29} :
\begin{eqnarray}
\beta_1(s)=\int_{-\infty}^{\infty} e^{i(s\xi+\xi^2)}\D\frac d{d\xi}
\left[e^{-i\frac{a^2}2\log |\xi|}\chi_{[0,+\infty)}(\xi)\right]d\xi,\label{eq1.30} \\\nonumber\\
\beta_2(s)=\int_{-\infty}^{\infty} e^{i(s\xi+\xi^2)}\D\frac d{d\xi}
\left[e^{-i\frac{a^2}2\log |\xi|}\chi_{[(-\infty,0]}(\xi)\right]d\xi.\label{eq1.31} 
\end{eqnarray}

Recall that $A^+_a=(A_{a,1},A_{a,2},A_{a,3})\,,\,A^-_a=(A_{a,1},-A_{a,2},-A_{a,3})$ so that
$$\cos \gamma=1-2A_{a,1}^2\quad;\quad \sin\D\frac{\gamma}2=A_{a,1}.$$
Then it suffices to compute $A_{a,1}=\lim\limits_{s\to\infty} T_1(s)$.

But
$$T_1(s)=1-|\theta_1|^2=1-|a_{11}\beta_1+a_{21}\beta_1|^2$$
with
$$a_{1,1}=a_{2,1}=\D\frac{a}{2\sqrt 2\beta_1^{\prime}(0)}.$$
Hence
$$1-T_1=\D\frac{a^2}{4\left(\beta_1^{\prime}(0)\right)^2}|\beta_1(s)+\beta_2(s)|^2.$$
But notice that
$$\lim_{s\to\infty}(\beta_1+\beta_2)(s)=
\int_{-\infty}^{\infty}e^{i\xi}\D\frac d{d\xi}\left[e^{-i\frac{a^2}2\log|\xi|}\right]d\xi$$
and
$$\beta_1^{\prime}(0)=\D\frac{a^2}2\int_0^{\infty}e^{i\xi^2-\frac{a^2}2\log|\xi|}d\xi.$$
After some calculations we get
\begin{equation}
\label{eq1.32}
A_{a,1}=e^{-\frac{\pi a^2}{2}}.
\end{equation}
Finally notice that the dispersive relation of $\beta_1$ and $\beta_2$ is
\begin{equation}
\label{eq1.33}
\phi(\xi)=\xi^2-\D\frac{a^2}2\log|\xi|,
\end{equation}
which indicates that the non-linearity only makes a logarithmic correction to the free evolution.

\subsection{Spirals}

\smallskip

In this section we will review some of the results obtained by S. Gutierrez and L. Vega in \cite{GV}. So far we have dealt with the I.V.P.
\begin{equation}
\left\{\begin{array}{ccc}
iu_t+u_{ss}+\D\frac{u}2\big(|u|^2+A(t)\big)&=&0\\ \\
u(s,0)&=&a\,\delta.
\end{array}\right.
\end{equation}
In this section we will be interested in
\begin{equation}
\label{eq2.2}
u(s,0)=a_1\,\mbox{p.v.}\,\D\frac 1s.
\end{equation}

On the other hand concerning the BF we have characterized all the possible solutions with respect to the scaling that preserve arc length. They were determined up to rotations by a free parameter, namely the curvature must be $\D\frac {a}{\sqrt t}$ and the torsion $\D\frac {s}{2t}$.

In order to capture the solutions of the BF which come from the initial condition \eqref{eq2.2} we have to change the ansatz. Assume $\mathcal{A}$ is a real antisymmetric matrix. Then consider solutions of the BF
\begin{equation}
\label{eq2.3}
\chi_t=\chi_s\,\land\,\chi_{ss}
\end{equation}
of the type
\begin{equation}
\label{eq2.4}
\chi(s,t)=e^{\mathcal{A}/2\log t}\,\sqrt t\,G\left(\frac{s}{\sqrt t}\right)\quad t>0.
\end{equation}
Hence $G$ has to be a solution of  ($e^{\mathcal{A}/2\log t}$ is a rotation)
\begin{equation}
\label{eq2.5}
(I+A)\,G-s\,G^{\prime}=2G^{\prime}\,\land\,G^{\prime\prime}\quad,\quad
|G^{\prime}|^2=1\,.\qquad \left(G^{\prime}\cdot G^{\prime\prime}=0\right).
\end{equation}
Multiplying on the left by $G^{\prime}$ we obtain
\begin{equation}
\label{eq2.6}
G^{\prime\prime}=\D\frac 12(I+\mathcal{A})\,G\,\land \,G^{\prime},
\end{equation}
which together with the conditions
\begin{eqnarray}
(I+A)\,G(0)\cdot G^{\prime}(0) & = & 0,\label{eq2.7} \\\nonumber\\
|G^{\prime}(0)|^2 & = & 1,\label{eq2.8} 
\end{eqnarray}
is equivalent to \eqref{eq2.5}. It is quite striaghforward to prove that \eqref{eq2.6} has global existence. Also it can be easily computed numerically. Again in this case the relevant feature is the behaviour of $G$ at $\pm$ infinity. This follows from the identity
$$\left(e^{-\mathcal{A}\log s}\,\D\frac{G(s)}s\right)_s=\D\frac
{(-\mathcal{A}+I)\,G-G^{\prime}}{s^2}=2\,\D\frac{c\,b}{s^2}.$$

Hence if we prove that the curvature is bounded we can proceed as in the case of the kink solutions. This is a consequence of the following elemental lemma. Consider without loss of generality
$$\mathcal{A}=\left(\begin{array}{ccc}  0&  -\mu &  0 \\ \mu & 0  & 0  \\ 0 &  0 & 0 \end{array}\right)\quad \mu\in\R.$$

\begin{lemma}\label{L1} If $G$ solves \eqref{eq2.5} then
\begin{equation}
\label{eq2.9}
|T^{\prime}|^2(s)=-\mu\, T_3(s)-\nu\quad\mbox{with}
\end{equation}
\begin{equation}
\label{eq2.10}
\nu=-\mu \,T_3(0)-\D\frac 14\left|(I+\mathcal{A})\,G(0)\right|^2.
\end{equation}
\end{lemma}

\underline{Proof} Differentiating in (\ref{eq2.5}) we get
$$\mathcal{A}\,T=s\,T^{\prime}+2\,T\,\land\,T^{\prime}.$$
Then
$$\begin{array}{c}
\mathcal{A}\,T\,\land\,T=-s\,T\,\land\,T^{\prime}-2\,(T\cdot T^{\prime})\,T+2\,T^{\prime\prime}\\\\
(\mathcal{A}T\,\land\,T)\cdot T^{\prime}=2\, T^{\prime\prime}\cdot T^{\prime}=\D\frac d{ds}\,|T^{\prime}|^2.
\end{array}$$

But $(\mathcal{A}\,T\,\land\,T)\cdot T^{\prime}=-\mu\, T^{\prime}_3$ which concludes the proof.

\bigskip

It turns out that if $a\neq 0$ the frame $\left(T,\D\frac{\mathcal{A}T}{|\mathcal{A}T|},
\D\frac{\mathcal{A}T\land T}{|\mathcal{A}T\,\land\,T|}\right)$ is the right one to study this family of solutions in the spirit of what we have seen for the kink solutions. This is a consequence of the identities:
\begin{eqnarray}
c(n-ib)=\left(\D\frac{dc^2}{ds}-2ic^2\left(\tau-\D\frac{s}{\tau}\right)\right)
\left(\D\frac{\mathcal{A}T\,\land\,T}{|\mathcal{A}T\,\land\,T|^2} -
i\D\frac{\mathcal{A}T}{|\mathcal{A}T|^2}\right).\label{eq2.11}
\end{eqnarray}
Let us consider $(y,h)$ given by
\begin{eqnarray}
 y:=\D\frac{dc^2}{ds}=\D\frac{d|T^{\prime}|^2}{ds}\label{eq2.12}\\
 \nonumber\\
 h:=c^2\left(\tau-\D\frac s2\right).\label{eq2.13}
\end{eqnarray}
The key point then is that $(y,h)$ solves the system
\begin{equation}
\label{eq2.14}
\left\{\begin{array}{l}
y^{\prime}=sh+g(|T^{\prime}|^2),\qquad
g(|T^{\prime}|^2)=2E(0)-\D\frac{(3|T^{\prime}|^2+\nu) (|T^{\prime}|^2+\nu)}{2}\\\\
h^{\prime}=-\D\frac s4 y,
\end{array}\right.
\end{equation}
with $\nu$  given in \eqref{eq2.10}.

The relation with NLS of the above comments is also clear. The ansatz \eqref{eq2.4} implies that $c$ and $\tau$ have to be selfsimilar i.e.
$$c(s,t)=\D\frac 1{\sqrt t}\,\widetilde{c}\left(\frac{s}{\sqrt t}\right)\quad;\quad
\tau(s,t)=\D\frac 1{\sqrt t}\,\widetilde{\tau}\left(\frac{s}{\sqrt t}\right)$$
Therefore we have to look for solutions of \eqref{eq3} with ${A}(t)=\D\frac{\nu}t$ of the type $\D\frac 1{\sqrt t}\,\widetilde{u}\left(\frac{s}{\sqrt t}\right)$. This gives us ODE
$$\widetilde{u}^{\prime\prime}-\D\frac i2\left(\widetilde{u}+s\,\widetilde{u}^{\prime}\right)+
\D\frac{\widetilde{u}}2\left(|\widetilde{u}|^2+\nu\right)=0.$$
Then we consider $\widetilde{u}(s)=f(s)e^{is^2/4}$ with
\begin{equation}
\label{eq2.15}
f^{\prime\prime}+i\D\frac s2 f^{\prime}+\D\frac f2\left(|f|^2+\nu\right)=0.
\end{equation}
The fundamental property of \eqref{eq2.15} is that it has a natural energy
\begin{equation}
\label{eq2.16}
E(s)=|f^{\prime}|^2+\D\frac 14\left(|f|^2+\nu\right)^2=E(0).
\end{equation}
Finally the way to analyze the solutions of (\ref{eq2.15}) is by defining $y$ and $h$ as
\begin{equation}
\label{eq2.17}
\D\frac y2+ih=\overline{f}f^{\prime}.
\end{equation}
It turns out that $(y,h)$ solves the system \eqref{eq2.14}.
The results can be found in \cite{GV}.

With respect to the I.V.P. with data \eqref{eq2.2} we have the following result.

\begin{proposition}\label{2.P.2} Let $\,\nu\ge 0\,$ and consider the I.V.P.
\begin{equation}
\label{eq2.18}
\left\{\begin{array}{rcl}
iu_t+u_{ss}+\D\frac{u}2\left(|u|^2+\D\frac{\nu}t\right)&=&0\qquad t>0\quad s\in\R,\\\\
u(x,0)&=&a_1\,\mbox{p.v.}\,\D\frac 1s\qquad a_1\in\C\backslash \{0\}.
\end{array}\right.
\end{equation}
Then either there is no weak solution $u$ for the I.V.P. \eqref{eq2.18} in the class
$$t^{1/2}u\,,\quad t^{3/2}|u|^2u\in L^{\infty}\left([0,+\infty),\mathcal{S}(\R)\right)$$
with $\,\lim\limits_{t\downarrow 0}u(s,t)=a_1 \,p.v.\left(\D\frac 1s\right)\,$ in $\,\mathcal{S}^{\prime}(\R)\,$ or there is more than one.
\end{proposition}

\section{Stability}

In this section we shall present the known results concerning the stability of the BF selfsimilar solutions $\chi_a$, which show that under some appropriate small regular perturbations the corner still remains at time zero. To this purpose we first study the perturbations of the filament function $u_a(s,t)=a\,\frac{e^{i\frac{s^2}{4t}}}{\sqrt{t}}$, solutions of 
\begin{equation}
\label{eqNLS}
\left\{\begin{array}{rcl}
iu_t+u_{ss}+\D\frac 12\left(|u|^2-\D\frac{|a|^2}{t}\right)u&=&0,\\
u(s,0)&=&a\,\delta.
\end{array}\right.
\end{equation}
The formation of a corner for the flow of curves $\chi_a$ corresponds to the limit  $a\,\delta$ of $u_a$ as time goes to zero.
The informations about the perturbations of $u_a$ will be gathered in the first subsection. Then in the last subsection we shall deal with the construction of new BF solutions. We shall define from a perturbation of $u_a$ a curvature and a torsion function. Then we shall construct a flow of curves, solution of the BF with these curvature and torsion function. Finally, we shall prove that the new flow of curves is close to the selfsimilar one, and in particular a singularity is still generated at time zero.

\subsection{Perturbations of NLS solutions with $\delta$-function as initial datum}

In an attempt to understand the above equation \eqref{eqNLS}, we considered in \cite{BaVe1} the equation 
\begin{equation}\label{dirac}
iu_t+\partial_{ss} u\pm |u|^{\alpha}u=0,
\end{equation}
with an initial data of Dirac type, and its solution
$$u(s,t)=u_{a,\pm\alpha}(s,t)=u_a(s,t)\,e^{\pm iA_{a,\alpha}(t)},$$
with
$$A_{a,\alpha}(t)=\left\{\begin{array}{rcl}
\frac{|a|^{\alpha}}{1-\frac{\alpha}{2}}t^{1-\frac{\alpha}{2}}\qquad \text{si}\quad \alpha\neq 2,\\
\\
|a|^{2}\, \log\,t\qquad \text{si}\quad \alpha=2.
\end{array}\right.$$
For $\,\alpha\geq 2$, $u_{a,\pm\alpha}$ does not have a limit as $t$ goes to zero, and this gives an instability result, as explained in the previous subsection 2.3.\\

A first natural question is to find out if $\,u_{a,\pm\alpha}\,$ is a stable solution in the case $\,\alpha<2$. Denoting
$$\eta(s,t)=e^{\mp iA_{\alpha,a}(t)}(u-u_{a,\pm\alpha})(s,t),$$ 
and using the fact that $u_a$ is a solution of the linear equation, we have
$$\left\{\begin{array}{ccc}
i\eta_t+\partial_{ss}\eta\pm\left(|\eta+u_a|^{\alpha}-|u_a|^{\alpha}\right)\left(\eta+u_a\right)=0,\\
\eta(s,0)=u_0(s).
\end{array}\right.$$
If $\alpha<2$ then the term $\,|u_a|^{\alpha}\,$ is integrable at $t=0$, and we can perform a fixed point argument. We have obtained the following positive answer, valid also in higher dimensions, for $\,\alpha<\frac 2d$ (a similar result is presented in \cite{BaVe1} for Sobolev spaces, with a more technical proof).

\begin{theorem} \label{thsubcrit} Let $\,\alpha<2$ and $\,u_0\in L^2(\mathbb R)$. There exists a time $\,T=T(a,\|u_0\|_{L^2})\,$ and a unique solution $\,u\,$ of equation \eqref{dirac} with initial data $u_0+a\,\delta$, such that
$$u-u_{a,\pm \alpha}\in L^4\left([0,T),L^\infty(\mathbb R)\right)\cap\mathcal{C}\left([0,T),L^2(\mathbb R)\right).$$
\end{theorem}

A case of rough perturbation of $a\delta_0$ has been studied by Kita. In \cite{Ki} he has described the solution generated by the sum of two or three Dirac masses.

Finally, we shall do a parallel with the corresponding parabolic problem,
$$\left\{\begin{array}{rcl}
u_t-\Delta u\pm |u|^{\alpha}u=0,\\
u(s,0)=a\,\delta,\end{array}\right.
$$
which has been intensively studied. F.B. Weissler has shown in \cite{Weissler} that in the focusing case, for $\alpha\leq\frac 2d$, there is no unicity of the solution. In the defocusing case, if only positive solutions are considered, H. Br\'ezis and A. Friedman have proven in \cite{BrFr} the existence of a unique solution for $\alpha<\frac 2d$, and non-existence for $\alpha\geq \frac 2d$. Therefore for the heat equation there is an important difference between the focusing and the defocusing case, even for $\alpha<\frac 2d$. Theorem \ref{thsubcrit} shows that this is not the case for the Schr\"odinger equation.\\

The second natural question is to find out if $u_a$ is a stable solution of the cubic 1-D Schr\"odinger equation modified by the complex phase $e^{\pm ia^2\log t}$.

Having as a purpose the study of perturbations of the selfsimilar flow $\chi_a$, we shall first consider perturbations at positive times of the filament function
$$u_a(s,t)=a\frac{e^{i\frac{s^2}{4t}}}{\sqrt{t}},$$
and we shall treat more generally the equation
$$iu_t+u_{ss}\pm\left(|u|^2-\frac{a^2}{t}\right)u=0$$
around its particular solution $u_a$. By defining $v$ the pseudo-conformal transformation of $u$,
\begin{equation}\label{calT}
u(s,t)=\D\frac{e^{i\frac{s^2}{4t}}}{\sqrt{t}}\overline{v}\left(\frac st,\frac 1t\right),
\end{equation}
this is equivalent to study the large time behavior of perturbations of the constant solution $a$ of 
\begin{equation}\label{GP}
iv_t+v_{ss}\pm \frac 1t\left(|v|^2-a^2\right)v=0.
\end{equation}
Let us notice that this equation has a natural energy
$$E(t)=\frac{1}{2}\int |v_s(t)|^2\,ds\mp\frac{1}{4t}\int(|v(t)|^2-a^2)^2\,ds.$$
with the law
$$\partial_{t}E(t)\mp\frac{1}{4t^2}\int(|v|^2-a^2)^2\,ds=0.$$
In the defocusing case, this has allowed us to obtain in \cite{BaVe1} a good control in time of the norm $\|v(t)-a\|_{L^2}$ and to conclude global existence.

As the coefficients of the linear terms are of the type $\frac 1t$, and since the perturbation is done around a constant function, there will be connections with the long range effects of 1-D cubic Schr\"odinger equation (T. Ozawa \cite{Ozawa}, R. Carles \cite{Carles}, N. Hayashi and P. Naumkin \cite{HaNa}) 
and with the scattering around the constant solution $1$ of 2-D Gross-Pitaevskii equation (S. Gustafson, K. Nakanishi and T.-P. Tsai \cite{GuNaTs}).

We shall search for wave operators for the perturbation $w(t)=v(t)-a$, meaning that we shall search for solutions of 
$$iw_t+w_{ss}\pm \frac 1t\left(|w+a|^2-a^2\right)(w+a)=0,$$
behaving at large times as a given profile $w_1(s,t)$. This can be done by a fixed point argument, in a suitable space defined around $w_1$, for the application
$$Aw=w_1\mp i\int_t^{\infty} e^{i(t-\tau)\partial_s^2}\left(\frac{(|w+a|^2-a^2)(w+a)}{\tau}-\frac{(|w_1+a|^2-a^2)(w_1+a)}{\tau}\right)d\tau$$
$$\mp i\int_t^{\infty} e^{i(t-\tau)\partial_s^2}\left(\frac{(|w_1+a|^2-a^2)(w_1+a)}{\tau}-(i\partial_\tau+\partial_x^2)w_1\right)d\tau.$$
For showing that $w$ behaves like a free evolution $e^{it\partial_s^2}\,u_+$, one should control the Duhamel terms corresponding to the worse source terms, the linear ones $\frac{a^2}{t}w_1$ and  $\frac{a^2}{t}\overline{w}_1$. These Duhamel terms are
$$a^2\int_t^\infty e^{i(t-\tau)\partial_s^2}\,e^{i\tau\partial_s^2}\,u_+\frac{d\tau}{\tau}\qquad,\qquad
 a^2\int_t^\infty e^{i(t-\tau)\partial_s^2}\,e^{-i\tau\partial_s^2}\,\overline{u}_+\frac{d\tau}{\tau}.$$
We notice that the first integral in not controllable. To avoid this term, we shall look for $w$ close to the long range ansatz 
$$w_1(s,t)=e^{\pm i a^2\log t}\,e^{it\partial_s^2}\,u_+(s).$$
By using inhomogeneus Strichartz inequalities we manage to control the remaining linear term
$$a^2\left\|\int_t^{\infty} e^{i(t-\tau)\partial_s^2}\left(e^{-i\tau\partial_s^2}\frac{\overline{u_+}}{\tau^{1\pm ia^2}}\right)d\tau\right\|_{L^2_s}=a^2\left\|\int_t^{\infty} e^{-i2\tau\partial_s^2}\left(\frac{\overline{u_+}}{\tau^{1\pm ia^2}}\right)d\tau\right\|_{L^2_s}$$
$$\leq Ca^2\,\left\|\frac{\overline{u_+}}{\tau}\right\|_{L^{{p'}}((t,\infty),L^{{q'}})}=Ca^2\,\|u_+\|_{L^{{q'}}}\left\|\frac{1}{\tau}\right\|_{L^{{p'}}(t,\infty)}=Ca^2\,\|u_+\|_{L^{{q'}}}\frac{1}{t^\frac{1}{p}}.$$ 
The 1-admissible couples $(p,q)$ are between $(\infty,2)$ and $(4,\infty)$, therefore if  $u_+\in L^1\cap L^2$, the best decay we can obtain this way is $t^{-\frac{1}{4}} $. We denote
$$v_1(s,t)=a+w_1(s,t)=a+e^{\pm i a^2\log t}\,e^{it\partial_s^2}\,u_+(s).$$
We get our first result of \cite{BaVe2}, in mixed norm spaces $L^pL^q$.
\begin{theorem}\label{longrange}
Let $t_0>0$. There exists a constant $a_0>0$ such that for all $a<a_0$, and for all $u_+$ small in $L^1\cap L^2$ with respect to $a_0$ and $t_0$, equation \eqref{GP} has a unique solution $v$ with 
$$v-v_1\in \mathcal{C}([t_0,\infty),L^2(\mathbb R))\cap L^4([t_0,\infty),L^\infty(\mathbb R)),$$ 
satisfying, as $t$ goes to infinity,
$$\|v(t)-v_1(t)\|_{L^2}+\|v-v_1\|_{L^4((t,\infty),L^\infty)}=\mathcal{O}(t^{-\frac 14}).$$
\end{theorem}

However, the $t^{-\frac 14}$ decay, optimal for the methods used in mixed norm spaces, is not enough to construct BF solutions. We shall then 
exploit more the oscillations of the linear remaining term, by using the Fourier transform and strengthening the conditions on $u_+$. Writting
$$\int_t^{\infty} e^{i(t-2\tau)\partial_s2}\left(e^{-i\tau\partial_s^2}\frac{\overline{u_+}}{\tau^{1\pm ia^2}}\right)d\tau=\int_t^{\infty} \int\frac{e^{-i(t-2\tau)\xi^2}}{\tau^{1\pm ia^2}}\,e^{ix\xi}\,\widehat{u_+}(-\xi)d\xi d\tau ,$$
and doing an integration by parts in $\tau$ we get
$$=\int\frac{e^{it\xi^2}}{t^{1\pm ia^2}} \,e^{is\xi}\,\frac{\widehat{u_+}(-\xi)}{2i\xi^2}d\xi+(1\pm ia^2)\int_t^{\infty}\int\frac{e^{-i(t-2\tau)\xi^2}}{\tau^{2\pm ia^2}} \,e^{is\xi}\,\frac{\widehat{u_+}(-\xi)}{2i\xi^2}d\xi d\tau.$$
Using Plancherel formula we obtain the control
$$\left\|a^2\int_t^{\infty} e^{i(t-\tau)\partial_s^2}\left(e^{-i\tau\partial_s^2}\frac{\overline{u_+}}{\tau^{1\pm ia^2}}\right)d\tau\right\|_{L^2}\lesssim \frac{a^2}{t}\left\|\frac{\widehat{u_+}(-\xi)}{\xi^2}\right\|_{L^2}=\frac{a^2}{t}\|u_+\|_{\dot{H}^{-2}}.$$
The derivatives $\nabla^k$ of this term have a $1/t$ decay, with constants depending on the $\dot{H}^{-2+k}$ norm of $u_+$. We obtain in \cite{BaVe2} a new result, with a better decay.

\begin{theorem}\label{longrangeH}
Let $t_0>0$ and $m\in\mathbb{N}^*$. There exists a constant $a_0>0$ such that for all $a<a_0$ and for all $u_+$ small in $\dot{H}^{-2}\cap H^{m}\cap W^{m,1}$ with respect to $a_0$ and $t_0$, equation (\ref{GP}) has a unique solution
$$v-v_1\in \mathcal{C}([t_0,\infty),H^m(\mathbb R)),$$
satisfying, as $t$ goes to infinity, for all entire number $0<k\leq m$,
\begin{equation}\label{rateH}\|(v-v_1)(t)\|_{L^2}=\mathcal{O}(t^{-\frac 12})\,\,\,,\,\,\,\|\partial_s^k(v-v_1)(t)\|_{L^2}=\mathcal{O}(t^{-1}).\end{equation}
\end{theorem}

The presence of two different rates of decay is now due to the quadratic terms $\frac{2a}{t}|w_1|^2$ and $\frac{a}{t}w_1^2$. The less oscillating is the first one, which gives the rate of decay $t^{-\frac 12}$. However, its derivatives can be better upper-bounded. We shall sketch in the following the way we treat the Duhamel term coming from the worse quadratic term , $\frac{|w_1|^2}{t}$. We shall use again the Fourier oscillations,
$$\partial_s^k\int_t^{\infty} e^{i(t-\tau)\partial_s^2}\left(\frac{\left|e^{i\tau\partial^2_s}u_+\right(\tau)|^2}{\tau}\right)d\tau=\int_t^{\infty}\int|\xi|^k\frac{e^{-i(t-\tau)\xi^2}}{\tau}\,e^{is\xi}\,\widehat{e^{i\tau\partial_s^2}u_+}* \widehat{e^{-i\tau\partial_s^2}\overline{u}_+}\,d\xi\, d\tau$$
$$=\int_t^{\infty} \int\int\frac{e^{-it\xi^2+2i\tau\xi(\eta-\xi)}}{\tau}\,|\xi|^k\,e^{is\xi}\,\widehat{u_+}(\eta) \widehat{\overline{u}_+}(\eta-\xi)\,d\eta\,d\xi \,d\tau,$$
and an integration by parts in $\tau$
$$=\int\int\frac{e^{-it(\xi^2+2\xi(\eta-\xi))}}{t}\,|\xi|^k\,e^{is\xi}\,\frac{\widehat{u_+}(\eta) \widehat{\overline{u}_+}(\eta-\xi)}{2i\xi(\eta-\xi)}\,d\eta\,d\xi$$
$$+\int_t^{\infty}\int\int\frac{e^{-it\xi^2-2i\tau\xi(\eta-\xi)}}{\tau^2}\,|\xi|^k\,e^{is\xi}\,\frac{\widehat{u_+}(\eta) \widehat{\overline{u}_+}(\eta-\xi)}{2i\xi(\eta-\xi)}\,d\eta\,d\xi\,d\tau.$$
Hence we have obtained that this quadratic Duhamel term can be controlled in $\dot{H}^1$ by
$$\frac{Ca}{t}\left\|\widehat{u_+}\star \frac{\widehat{\overline{u_+}}}{\cdot}\right\|_{L^2}=\frac{Ca}{t}\left\| u_+ \,\mathcal F\left(\frac{\widehat{\overline{u_+}}}{\cdot}\right)\right\|_{L^2}\leq \frac{Ca}{t}\|u_+\|_{H^1}\|u_+\|_{\dot{H}^{-1}}.$$
We notice that this argument can be repeated for $\dot{H}^k$, with $k>1$, but not for $k=0$.

At the level of the initial equation before the pseudo-conformal transformation, we obtain that $u(t)$ is close to
$$u_1(s,t)=u_a+e^{\pm ia^2\log t}\,\hat {\overline{u}}_+\left(-\frac s2\right).$$
The condition $u_+\in\dot{H}^{-2}$ is translated on $\hat {\overline{u}}_+$ by a canceling condition at the origin, which could be related to the blow-up result of J. Bourgain and W. Wang \cite{BoWa}. We notice the absence of a limit in $L^2$ for $u(t)-u_a$ at time zero. However this does not exclude the possibility to construct from $u$ a binormal flow with a trace at zero. This is because the curvature and the torsion defined from $u$ corresponds to the second derivatives of the curves.

\subsection{On the stability of selfsimilar solutions of the binormal flow}
  
In this last subsection we sketch how from a solution of Theorem \ref{longrangeH} we construct a flow of curves $\chi$, solution of the BF, and get information about its limit as $t$ goes to zero. We prove that the new curve $\chi$ is close to the selfsimilar one $\chi_a$, and then obtain the persistence of the formation of a singularity at time zero.

Let $v$ be a solution of Theorem \ref{longrangeH} with $m=3$, and $u$ its pseudo-conformal transform \eqref{calT}. Since $v$ is regular enough and does not vanish, we can define two real functions $c$ and $\phi$ such that 
\begin{equation}\label{def}
u(s,t)=c(s,t)e^{i\phi(s,t)}.
\end{equation}
We define $\tau(s,t):=\phi_s(s,t)$, so
$$u(s,t)=c(s,t)e^{i\int_0^s\tau(s,t)ds +\phi(0,t)}.$$
Then, the function
$$\tilde{u}(s,t)=c(s,t)e^{i\int_0^s\tau(s,t)ds}$$
is a filament function, solution of \eqref{eq3} with $A(t)$ replaced by $\frac{a^2}{t}+\phi_t(0,t)$. In particular, $(c,\tau)$ satisfy the intrinsic equations \eqref{eq1.41}, and a computation yields the extra-information at $s=0$
\begin{equation}\label{extra}
\frac{a^2}{t}+\partial_t\phi(0,t)=2\left(\frac{c_{ss}-c\tau^2}{c}\right)(0,t)+c^2(0,t).
\end{equation}
In view of the pseudo-conformal transform, the curvature is
 $$c(s,t)=|u(s,t)|=\frac{1}{\sqrt{t}}\left|v\left(\frac st,\frac 1t\right)\right|,$$
and the torsion is
$$\tau(s,t)=\Im\frac{u_s(s,t)}{u(s,t)}=\Im\frac{\frac{is}{2t}\overline{v}\left(\frac st,\frac 1t\right)+\partial_s \overline{v}\left(\frac st,\frac 1t\right)}{\overline{v}\left(\frac st,\frac 1t\right)}.$$
As $v\in H^3(\mathbb R)\subset\mathcal{C}^\frac 52(\mathbb R)$, the curvature and torsion function will be regular enough. 

By definition,
$$v(s,t)=a+e^{\pm i a^2\log t}e^{it\partial_s^2}u_+(x)+(v-v_1)(s,t).$$
The decay of Theorem \ref{longrangeH} combined with the explicit formula of the free Schr\"odinger evolution allows us to obtain that $(c,\tau)$ is close to $(c_a,\tau_a)$, with strong rate of decay as $t$ goes to zero, and constants depending on the norms of $u_+$.

Starting from $(c,\tau)$ we construct the binormal flow as follows. We fix an initial data $(T,n,b)(0,\tilde{t}_0)$ at time $\tilde{t}_0=1/t_0$, to be specified later. For $\tilde{t}_0>t>0$ we define $(T,n,b)(0,t)$ by imposing that the derivatives in time satisfy the same system as the tangent, normal and binormal vectors of a binormal flow of curves satisfy in general. Then, we define $(T,n,b)(s,t)$ from $(T,n,b)(0,t)$ by integrating the Frenet system for fixed time $t$. We therefore impose
$$\left(\begin{array}{c} T\\ n\\ b \end{array}\right)(s,t)=
\left(\begin{array}{c} T\\ n\\ b \end{array}\right)(0,\tilde{t}_0)-
\int_t^{\tilde{t}_0}
\left(\begin{array}{ccc} 0 & -c\,\tau & c_x \\ c\,\tau & 0 & \left(\frac{c_{ss}-c\tau^2}{c}\right) \\ -c_s &  -\left(\frac{c_{ss}-c\tau^2}{c}\right) & 0 \end{array}\right)
\left(\begin{array}{c} T\\ n\\ b \end{array}\right)(0,t')\,dt'$$
$$+
\int_0^s\left(\begin{array}{ccc} 0 & c & 0 \\ -c & 0 & \tau \\ 0 &  -\tau & 0 \end{array}\right)
\left(\begin{array}{c} T\\n\\b \end{array}\right)(s',t)\,ds'.$$
Proceeding this way, $T$ shall satisfy the equation
$$T_t=T\land T_{ss},$$
of the tangent vector of a general binormal flow. 

Once $T$ is constructed, for a given $\chi(0,\tilde{t}_0)$ we define
$$\chi(s,t)=\chi(0,\tilde{t}_0)-\int^{\tilde{t}_0}_t (c\,b)(0,t')\,dt'+\int_0^sT(s,t)\,ds.$$
By using Frenet system, 
$$T_t=T\land T_{ss}=T\land (c\,n)_s=T\land(c_s\, n +c\,\tau\, b)=-c\,\tau\, n+c_s\,b,$$
and it follows that $\chi$ satisfy the binormal flow. 

Finally, $\chi(s,t)$ is obtained this way for all times when $(c,\tau)$ are regular, that is $\tilde{t}_0>t>0$. Knowing that $c$ is close to $c_a=\frac{a}{\sqrt{t}}$, we get
$$\left|\chi(s,t_1)-\chi(s,t_2)\right|=
\left|\int_{t_1}^{t_2} c(s,t)\,b(s,t)\, ds\right|\leq \int_{t_1}^{t_2} \frac{Ca}{\sqrt{t}}\, dt\underset{t_1,t_2\rightarrow 0}{\longrightarrow} 0.$$
This implies the existence of a limit $\chi_0(s)$ at $t=0$, and we obtain similarly that for all $s\in(-\infty,\infty)$,
$$\left|\chi(s,t)-\chi_0(s)\right|\leq Ca\sqrt{t} .$$

Next, for proving that $\chi(s,0)$ is close to $\chi_a(s,0)$, we first show that the tangent $T(s,t)$ is close to $T_a(s,t)$. As a first step, we prove that $T(0,t)$ is close to $T_a(0,t)$, and then we prove that $T(s,t)$ is close to $T_a(s,t)$. For the second step we need to know that $(T,n,b)(0,t)$ is close to $(T_a,n_a,b_a)(0,t)$. In the following we shall sketch the estimates at $(0,t)$, which is the less technical step, but nevertheless a key point of the proof.

From the above definition of $(T,n,b)$, 
$$\left(\begin{array}{c} T\\ n\\ b \end{array}\right)_t(0,t)=
\left(\begin{array}{ccc} 0 & -c\,\tau & c_s \\ c\,\tau & 0 & \left(\frac{c_{ss}-c\tau^2}{c}\right) \\ -c_s &  -\left(\frac{c_{ss}-c\tau^2}{c}\right) & 0 \end{array}\right)
\left(\begin{array}{c} T\\ n\\ b \end{array}\right)(0,t).$$
Now we use the extra-information at $x=0$ \eqref{extra}, 
$$\frac{c_{ss}-c\tau^2}{c}=\frac{c_a^2-c^2}{2}+\frac{\phi_t}{2}.$$
Therefore we obtain
$$\left(\begin{array}{c}
T\\n\\b
\end{array}\right)_t(0,t)=
\left(\begin{array}{ccc}
0 & -c\,\tau & c_s \\ c\,\tau & 0 &  \frac{c_a^2-c^2}{2}+\frac{\phi_{t}}{2}
\\ -c_s &  -\frac{c_a^2-c^2}{2}-\frac{\phi_{t}}{2} & 0 
\end{array}\right)
\left(\begin{array}{c}
T\\ n\\ b
\end{array}\right)(0,t).$$
In order to get rid of the term $\phi_t(t,0)$, we perform the change of functions $\tilde{n}+i\tilde{b}=e^{i\frac{\phi}{2}}(n+ib)$, and we impose as an initial condition $(T,\tilde{n},\tilde{b})(0,\tilde{t}_0)=(T_a,n_a,b_a)(0,\tilde{t}_0)$. For having $(T,\tilde{n},\tilde{b})(0,t)$ close to $(T_a,n_a,b_a)(0,t)$, it is then enough to show that the matrix
$$\left(\begin{array}{ccc} 0 & |c\,\tau| & |c_s| \\ |c\,\tau| & 0 & \left|\frac{c_a^2-c^2}{2}\right| \\ |c_s| & \left|\frac{c_a^2-c^2}{2}\right| & 0 \end{array}\right)(t,0)$$
is integrable in time near $t=0$. This is true, in view of the strong rate of decay obtained for $(c,\tau)-(c_a,\tau_a)$. Also, we get a good decay for $|e^{i\frac{\phi}{2}}-1|$, so $(T,\tilde{n},\tilde{b})(0,t)$ is close to $(T,n,b)(0,t)$. In conclusion, we have $(T,n,b)(0,t)$ close to $(T_a,n_a,b_a)(0,t)$.\\

We obtain finally that for $\epsilon>0$, and $u_+$ small enough, 
$$|(T-T_a)(s,t)|\leq\epsilon.$$
In particular, for $s>0$,
$$\chi(s,t)-\chi(0,t)=\int_0^s T (s,t) \,ds=A^+_a s +\int_0^s (T-T_a)(s,t)\,ds +\int_0^s T_a(s,t)-A^+_a \,ds.$$
The third term tends to zero as $t$ goes to zero, uniformly in $s$ \cite{GuRV}. By letting $t$ go to zero, we obtain the desired information
$$|\chi_0(s)-\chi_0(0)-sA^+_a|\leq \epsilon\, s.$$ 

More precisely, we have obtained in \cite{BaVe2} the following result.
\begin{theorem}\label{stab}
Let $\epsilon>0$, $t_0>0$, $0<a<a_0$, where $a_0$ is the constant in Theorem \ref{longrangeH}. Let $u_+$ small in $\dot{H}^{-2}\cap H^3\cap W^{3,1}$ with $s^2u_+$ small in $H^1$ with respect to $\epsilon$, $a$ and $t_0$, and let $v$ be the corresponding solution obtained in Theorem \ref{longrangeH}. By using the Hasimoto transform, we construct from $v$ a family of curves $\chi(s,t)$ which solves the binormal flow for $\tilde{t_0}>t>0$, and such that there exists a unique $\chi_0$ satisfying
$$|\chi(s,t)-\chi_0(s)|<Ca\sqrt t$$
uniformly in $s\in(-\infty,\infty)$.\\
Moreover,
$$|\chi_0(s)-\chi_0(0)-sA^+_a \mathbb{I}_{[0,\infty)}(s)-sA^-_a \mathbb{I}_{(-\infty,0]}(s)|<\epsilon \,|s|,$$
where $A^\pm_a$ are two distinct non opposite unit vectors, such that the angle $\gamma$ between $A_a^+$ and $-A_a^-$ is determined by the relation 
$$\sin\frac\gamma 2=e^{-\frac{\pi a^2}{2}}.$$ 
\end{theorem}
Since
$$\chi_a(s,0)=sA^+_a \mathbb{I}_{[0,\infty)}(s)+sA^-_a \mathbb{I}_{(-\infty,0]}(s),$$
we have that $\chi_0(s)$ lies in the $\epsilon$-cone around $\chi_a(s,0)$. Therefore we have obtained a large family of curves, solutions of the binormal flow, which are close to the selfsimilar solution $\chi_a$. In particular a corner is still formed for $\chi(s,t)$ at $t=0$ and at $s=0$. The angle of this corner can be made as close as desired to the one in between $A^+_a$ and $-A^-_a$  by taking $\epsilon$ small enough. 
This shows that the formation of a singularity in finite time of $\chi_a$ is not an isolated phenomenon. In a work in progress we study, instead of the existence of wave operators in Theorem \ref{longrangeH}, the asymptotic completeness question, which could lead to a better description of the allowed perturbations.

\end{document}